# Directional Differentiability of the Generalized Metric Projection in Banach Spaces

Jinlu Li

Department of Mathematics
Shawnee State University
Portsmouth Ohio, 45662 USA

**Abstract** Let $X$ be a real uniformly convex and uniformly smooth Banach space and $C$ a nonempty closed and convex subset of $X$. Let $\Pi_C: X \to C$ denote the generalized metric projection operator introduced by Alber in [1]. In this paper, we define the Gâteaux directional differentiability of $\Pi_C$. We investigate some properties of the directional differentiability of $\Pi_C$. In particular, if $C$ is a closed ball, or a closed and convex cone (including proper closed subspaces), or a closed and convex cylinder, then, we give the exact representations of the directional derivatives of $\Pi_C$. We also compare the differences of the directional differentiability between $\Pi_C$ and the (standard) metric projection operator $P_C$.

## 1 Introduction

Let $(X, \|\cdot\|)$ be a real uniformly convex and uniformly smooth Banach space with topological dual space $(X^*, \|\cdot\|_*)$. Let $C$ be a nonempty closed and convex subset of $X$. Let $P_C: X \to C$ denote the (standard) metric projection operator. For any $x \in X$, $P_C x \in C$ such that

$$\|x - P_C x\| \leq \|x - z\|, \text{ for all } z \in C. \tag{1.1}$$

For uniformly convex and uniformly smooth Banach space $X$, $P_C: X \to C$ is a well-defined single-valued mapping. The point $P_C x \in C$ is considered as the best approximation of $x$ by elements of $C$, which is the closest point from $x$ to $C$.

Let $J: X \to X^*$ be the normalized duality mapping in uniformly convex and uniformly smooth Banach space $X$, which is a well-defined single-valued mapping. With the help of $J$, the basic variational principle of $P_C$ is: for any $x \in X$ and $y \in C$,

$$y = P_C x \quad \Leftrightarrow \quad \langle J(x - y), y - z\rangle \geq 0, \text{ for all } z \in C. \tag{1.2}$$

In particular, if $X$ is a Hilbert space, then normalized duality mapping $J$ coincides with the identity map in $X$ and the basic variational principle of $P_C$ becomes: for any $x \in X$ and $y \in C$,

$$y = P_C x \quad \Leftrightarrow \quad \langle x - y, y - z\rangle \geq 0, \text{ for all } z \in C. \tag{1.3}$$

In 1996, Alber in [1] introduced the concepts of generalized metric projection $\Pi$ by using a Lyapunov functional $V: X \times X \to \mathbb{R}_+$ in uniformly convex and uniformly smooth Banach spaces (see (2.3)). He proved that in uniformly convex and uniformly smooth Banach spaces, the generalized metric projection $\Pi_C: X \to C$ is a well-defined single-valued mapping that satisfies the following basic variational principle: for any $x \in X$ and $y \in C$,

$$y = \Pi_C(x) \iff \langle Jx - Jy,\ y - z\rangle \geq 0, \text{ for all } z \in C. \tag{1.4}$$

In [14], the present author extended the concepts of both the generalized projection and the generalized metric projection from uniformly convex and uniformly smooth Banach spaces to reflexive Banach spaces. In [2], Alber and Li showed that, even in uniformly convex and uniformly smooth Banach spaces, the generalized metric projection is different from the metric projection, in general. Recently, in [11], Khan, Li and Reich studied the metric projection, generalized projection and the generalized metric projection in general Banach spaces. In this paper, the authors showed that there are a great differences between the metric projection and the generalized metric projection in general Banach spaces.

Both the (standard) metric projection $P$ and the generalized metric projection $\Pi$ have been widely applied to approximating theory, optimization theory, theory of variational inequalities, fixed point theory, and so forth. In the applications, the basic variational principles play very important roles respectively in Hilbert spaces and in Banach spaces.

There are two important issues regarding to the applications of the standard metric projection $P$ and the generalized metric projection $\Pi$. The first issue is about their continuity, which has been studied by many authors (see [1, 5, 8, 22]). The second issue is a more complicated issue that is about the directional differentiability. The basic variational principle and the nonexpansive property of $P_C$ in Hilbert spaces have been used to define and to study the G$\hat{a}$teaux type directional differentiability of $P_C$ in Hilbert spaces (see [6−10, 16, 19]). In uniformly convex and uniformly smooth Banach spaces, the basic variational principle of the metric projection $P$ closely depends on the normalized duality mapping $J$, which is not a linear operator, in general. So, it is more complicated to study the G$\hat{a}$teaux directional differentiability of $P_C$ in Banach spaces (see [2, 4, 6, 17, 18, 21, 23]). Furthermore, some applications of the G$\hat{a}$teaux directional differentiability of $P_C$ has been applied to approximation problems, convex programming, theory of optimal control, and so forth (see [3, 17, 18, 20]).

Very recently, in [13] the present author studied the G$\hat{a}$teaux directional differentiability of $P_C$ in uniformly convex and uniformly smooth Banach spaces, which followed the steps of differentiability of functions in calculus. At first, the definition of the G$\hat{a}$teaux directional differentiability of $P_C$ is given. By this definition, some properties of the G$\hat{a}$teaux directional derivatives $P_C'$ of $P_C$ are proved. Then, for some special cases of the considered closed and convex subset $C$, such as closed balls, closed and convex cones, the specific formulas (solutions) for $P_C'$ are provided. Following the work in [13], in [15], Li, Cheng and others concentrated in studying the properties and solutions of $P_C'$ for Hilbert spaces and Hilbert-Bochner spaces. In [16] the present author studied the properties and solutions of $P_C'$ in uniformly convex and uniformly smooth Bochner spaces.

So far, as far as I know, the directional differentiability of the generalized metric projection operator $\Pi_C$ has never been studied by any author. In this paper, we deal with the G$\hat{a}$teaux directional differentiability of the generalized metric projection $\Pi$ in uniformly convex and uniformly smooth Banach spaces. Similar to the results of the G$\hat{a}$teaux directional differentiability of the metric projection $P$, in this paper, we will give the definition of the G$\hat{a}$teaux directional differentiability of $\Pi$. Then, we investigate the properties of the G$\hat{a}$teaux directional differentiability of $\Pi$ and find some specific formulas for $\Pi_C'$ when $C$ is a closed ball, or, a closed and convex cone, or a closed and convex cylinder.

For any $y \in C$, similar to the inverse image $P_C^{-1}(y)$ of $y$ by the metric projection $P_C$ in $X$, we define the inverse image of $y$ by the generalized metric projection $\Pi_C$ in $X$ by

$$\Pi_C^{-1}(y) = \{x \in X\text{: } \Pi_C x = y\}. \tag{1.5}$$

We will find that, in contrast with the inverse image $P_C^{-1}(y)$ of $y$, the inverse image $\Pi_C^{-1}(y)$ of $y$ is more complicated. Then, it is easy to imagine that the solutions for $\Pi_C'$ must be much more complicated than the solutions for $P_C'$. In section 6, we concentrate to study the properties of $\Pi_C'$ in the real uniformly convex and uniformly smooth Banach spaces $l_p$, for some $1 < p < \infty$, from which, we show the connections between $P_C'$ and $\Pi_C'$.

## 2 Preliminaries

### 2.1 The function of smoothness of smooth Banach spaces

The important theme of this paper is about the directional differentiability of the generalized metric projection in uniformly convex and uniformly smooth Banach spaces. It describes the smoothness of the generalized metric projection, which must relate to the smoothness of the considered Banach spaces. It leads us to introduce the definition of function of smoothness of smooth Banach spaces, which will be multiply used in this paper.

Let $(X, \|\cdot\|)$ be a uniformly convex and uniformly smooth Banach space with topological dual space $(X^*, \|\cdot\|_*)$ (In this paper, all considered Banach spaces are real). Let $\langle\cdot, \cdot\rangle$ denote the real canonical evaluation pairing between $X^*$ and $X$. Let $S(X)$ denote the unit sphere of $X$. It is well known that $X$ is uniformly smooth, if and only if, the following limit

$$\lim_{t\downarrow 0}\frac{\|x+tv\|-\|x\|}{t},$$

exists uniformly for all $(x, v) \in S(X)\times S(X)$. Then, in [13], we introduced the following definition of smoothness of uniformly smooth Banach spaces.

**Definition 2.1.** Let $X$ be a uniformly convex and uniformly smooth Banach space. Define $\Psi$: $X\times(X\backslash\{\theta\}) \to \mathbb{R}$ by

$$\Psi(x, v) = \lim_{t\downarrow 0}\frac{\|x+tv\|-\|x\|}{t}, \text{ for any } (x, v) \in X\times(X\backslash\{\theta\}). \tag{2.1}$$

$\Psi_X$ is called the (extended) function of smoothness of this uniformly smooth Banach space $X$.

**Lemma 2.2.** *Let X be a uniformly convex and uniformly smooth Banach space. Then, for any* $(x, v) \in X\times(X\backslash\{\theta\})$, *we have*

$$\Psi(x, v) = \begin{cases} \|v\|, & \text{for } x = \theta, \\ \|v\|\Psi\left(\frac{x}{\|x\|}, \frac{v}{\|v\|}\right), & \text{for } x \neq \theta. \end{cases} \tag{2.2}$$

*In particular,*

$$\Psi(x, x) = \|x\|, \text{ for any } x \neq \theta.$$

The normalized duality mapping $J$ in Banach spaces plays a very important role in the analysis of Banach spaces. The normalized duality mapping $J$ on a uniformly convex and uniformly smooth Banach space $X$ with dual space $X^*$ is a single-valued mapping from $X$ to $X^*$, which is defined by $Jx \in X^*$ such that

$$\langle Jx, x\rangle = \|Jx\|_*\|x\| = \|x\|^2 = \|Jx\|_*^2, \text{ for any } x \in X.$$

We list some notations below used in this paper. For any $u, v \in X$ with $u \neq v$, we write

(a) $\overline{v,u}$ = {$tv$ + (1−$t$)$u$: $0 \leq t \leq 1$}, a closed segment in $X$ with ending points at $u$ and $v$;
(b) $\overline{v,\vec{u}}$ = {$v + t(u - v)$: $0 \leq t < \infty$}, a closed ray in $X$ with ending point at $v$;
(c) $\overleftrightarrow{v,u}$ = {$v + t(u - v)$: $-\infty < t < \infty$}, a line in $X$ passing through points $v$ and $u$.

**Corollary 2.8 in [12]** *Let X be a uniformly convex and uniformly smooth Banach space and K a closed cone in X with vertex at $v \in X$. Then*

(i) *If $v = \theta$, we have*

(a) *JK is a closed cone in X* with vertex at $\theta^* = J\theta$;*
(b) *However, in general*

$$K \text{ is convex } \nRightarrow JK \text{ is convex};$$

(ii) *If $v \neq \theta$ or K is a ray with $\theta \notin K$, then JK is not a cone (not a ray).*

## 2.2. The generalized metric projection Π in uniformly convex and uniformly smooth Banach Spaces

The concepts of generalized projection and generalized metric projection were introduced by Alber in [1] in uniformly convex and uniformly smooth Banach Spaces, which have been extended to general Banach spaces. Let $X$ be a uniformly convex and uniformly smooth Banach space with dual space $X^*$. A Lyapunov functional $V$: $X \times X \to \mathbb{R}_+$ is defined by the following formula:

$$V(x, y) = \|Jx\|_{X^*}^2 - 2\langle Jx, y\rangle + \|y\|_X^2$$
$$= \|x\|_X^2 - 2\langle Jx, y\rangle + \|y\|_X^2, \text{ for any } x, y \in X.$$

Let $C$ be a nonempty closed and convex subset of $X$. The generalized metric projection $\Pi_C$: $X \to C$ is a well-defined single-valued mapping that satisfies

$$V(x, \Pi_C x) = \inf_{y \in C} V(x, y), \text{ for any } x \in X. \tag{2.3}$$

That is, $\Pi_C x$ satisfies

$$\|Jx\|_{X^*}^2 - 2\langle Jx, \Pi_C x\rangle + \|\Pi_C x\|_X^2 = \inf_{y \in C}(\|Jx\|_{X^*}^2 - 2\langle Jx, y\rangle + \|y\|_X^2),$$

or,

$$\|x\|_X^2 - 2\langle Jx, \Pi_C x\rangle + \|\Pi_C x\|_X^2 = \inf_{y \in C}(\|x\|_X^2 - 2\langle Jx, y\rangle + \|y\|_X^2). \tag{2.4}$$

The generalized metric projection $\Pi_C$ has many properties (see [1]). The following variational principle of $\Pi_C$ will be used in this paper: For any given $x \in X$ and $y \in C$, we have

$$y = \Pi_C(x) \iff \langle Jx - J(y),\ y - z\rangle \geq 0, \text{ for all } z \in C. \tag{1.3}$$

The metric projection $P_C$ and the generalized metric projection $\Pi_C$ in a uniformly convex and uniformly smooth Banach space $X$ have the following connections.

(a) In general, $\Pi_C \neq P_C$;
(b) When $X$ is a Hilbert space, then $\Pi_C = P_C$.

By the definition of the inverse image of points by the generalized metric projection $\Pi_C$ in $X$ recalled in (1.5), we have the following results in [12].

**Theorem 3.5 in [12]**. *For $y \in C$, suppose that there is $x \in X$ with $x \neq y$ such that $y = \Pi_C(x)$. Then,*

(a) $\Pi_C^{-1}(y)$ *is a closed subset in* $X$;
(b) *In general,* $\Pi_C^{-1}(y)$ *is not convex*;
(c) *In general,* $\Pi_C^{-1}(y)$ *is not a cone.*

## 3. Directional differentiability of the generalized metric projection Π onto closed and convex subsets

In this paper, we always, otherwise stated, let $X$ be a uniformly convex and uniformly smooth Banach space with dual space $X^*$ and $C$ a nonempty closed and convex subset of $X$. The generalized metric projection $\Pi_C$ is a well-defined single valued mapping from $X$ onto $C$. Similar to the definition of the Gâteaux directional differentiability of the metric projection $P_C$ given in [13], we define the Gâteaux directional differentiability of $\Pi_C$ in the space $X$.

**Definition 3.1** For $x \in X$ and $v \in X$ with $v \neq \theta$, if the following limit exists (that is a point in $X$),

$$\lim_{t \downarrow 0} \frac{\Pi_C(x+tv)-\Pi_C(x)}{t}$$

then, $\Pi_C$ is said to be (Gâteaux) directionally differentiable at point $x$ along direction $v$, which is denoted by

$$\Pi_C'(x; v) = \lim_{t \downarrow 0} \frac{\Pi_C(x+tv)-\Pi_C(x)}{t}.$$

$\Pi_C'(x; v)$ is called the (Gâteaux) directional derivative of $\Pi_C$ at point $x$ along direction $v$; and $v$ is called a (Gâteaux) differentiable direction of $\Pi_C$ at $x$. If $\Pi_C$ is (Gâteaux) directionally differentiable at point $x \in X$ along every direction $v \in X$ with $v \neq \theta$, then $\Pi_C$ is said to be (Gâteaux) directionally differentiable at this point $x \in X$. It is denoted by

$$\Pi_C'(x)(v) = \lim_{t \downarrow 0} \frac{\Pi_C(x+tv)-\Pi_C(x)}{t}, \text{ for } v \in X \text{ with } v \neq \theta.$$

$\Pi_C'(x)(v)$ is called the (Gâteaux) directional derivative of $\Pi_C$ at point $x$ along direction $v$. Let $D$ be a nonempty subset in $X$. If $\Pi_C$ is (Gâteaux) directionally differentiable at every point $x \in D$, then $\Pi_C$ is said to be (Gâteaux) directionally differentiable on $D \subseteq X$.

By Definition 3.1, we give some properties of the directional differentiability of the generalized metric projection $\Pi_C$ onto nonempty closed and convex subsets of $X$. We see that these properties of $\Pi_C'$ below are very similar to the properties of $P_C'$, which are proved in [13].

**Lemma 3.2** *Let $u, w \in C$ with $u \neq w$. Then, $\Pi_C$ is directionally differentiable at $u$ and $w$ along directions $w - u$ and $u - w$, respectively, such that*

$$\Pi_C'(u; w-u) = w - u \quad \text{and} \quad \Pi_C'(w; u-w) = u - w.$$

*Proof*. Since $C$ is convex, then, for any positive number $\alpha$ with $0 < \alpha < 1$, we have $u + \alpha(w-u) \in \overline{u,w} \subseteq C$. We calculate

$$\Pi_C'(u; w-u)$$

$$= \lim_{t\downarrow 0} \frac{\Pi_C(u+t(w-u))-\Pi_C(u)}{t}$$
$$= \lim_{t\downarrow 0, t<1} \frac{u+t(w-u)-u}{t}$$
$$= \lim_{t\downarrow 0, t<1} \frac{t(w-u)}{t}$$
$$= \lim_{t\downarrow 0, t<1} \frac{t(w-u)}{t}$$
$$= w - u.$$

Rest of the proofs are similar to the above proof. □

**Lemma 3.3** *The following statements are equivalent*

(i) $\Pi_C$ *is directionally differentiable on X such that, for every point* $x \in X$,

$$\Pi'_C(x)(v) = \theta, \text{ for any } v \in X \text{ with } v \neq \theta;$$

(ii) $\Pi_C$ *is a constant operator*; *that is, C is a singleton.*

*Proof*. By Lemma 3.2, the proof of this lemma is similar to the proof of Lemma 4.4 in [13] and it is omitted here. □

Next lemma proves that the directional derivative of $\Pi_C$ is positive homogenous.

**Lemma 3.4** *For* $x \in X$ *and* $v \in X$ *with* $v \neq \theta$, *if* $\Pi_C$ *is directionally differentiable at x along direction v*, *then*, *for any* $\lambda > 0$, $\Pi_C$ *is directionally differentiable at x along direction* $\lambda v$ *satisfying*

$$\Pi'_C(x; \lambda v) = \lambda \Pi'_C(x; v), \text{ for any } \lambda > 0.$$

*Proof*. For $x \in X$ and $v \in X$ with $v \neq \theta$, suppose that $\Pi'_C(x; v) = \lim_{t\downarrow 0} \frac{\Pi_C(x+tv)-\Pi_C(x)}{t}$ exists. Then, for any $\lambda > 0$, we have

$$\Pi'_C(x; \lambda v) = \lim_{t\downarrow 0} \frac{\Pi_C(x+t\lambda v)-\Pi_C(x)}{t}$$
$$= \lambda \lim_{t\downarrow 0} \frac{\Pi_C(x+t\lambda v)-\Pi_C(x)}{\lambda t}$$
$$= \lambda \Pi'_C(x; v).$$ □

**Proposition 3.5** *Let* $y \in C$. *Suppose* $(\Pi_C^{-1}(y))^{\circ} \neq \emptyset$. *Then*, $\Pi_C$ *is directionally differentiable on* $(\Pi_C^{-1}(y))^{\circ}$ *such that*, *for any* $x \in (\Pi_C^{-1}(y))^{\circ}$, *we have*

$$\Pi'_C(x)(v) = \theta, \text{ for every } v \in X \text{ with } v \neq \theta.$$

*Proof*. For any given point $x \in (\Pi_C^{-1}(y))^{\circ}$, there is $\delta_x > 0$ such that

$$x + w \in (\Pi_C^{-1}(y))^{\circ}, \text{ for any } w \in X \text{ with } \|w\| < \delta_x.$$

That is,

$$\Pi_C(x + w) = \Pi_C(x) = y, \text{ for any } w \in X \text{ with } \|w\| < \delta_x. \tag{3.1}$$

For any $v \in X$ with $v \neq \theta$, by (3.1), we calculate

$$\begin{aligned}\Pi'_C(x)(v) &= \lim_{t\downarrow 0}\frac{\Pi_C(x+tv)-\Pi_C(x)}{t}\\ &= \lim_{t\downarrow 0,\, t<\delta_x\backslash\|v\|}\frac{\Pi_C(x+tv)-y}{t}\\ &= \lim_{t\downarrow 0,\, t<\delta_x\backslash\|v\|}\frac{\Pi_C(x)-y}{t}\\ &= \lim_{t\downarrow 0,\, t<\delta_x\backslash\|v\|}\frac{y-y}{t}\\ &= \theta.\end{aligned}$$

□

**Proposition 3.6** *Suppose* $C^{o} \neq \emptyset$. *Then* $\Pi_C$ *is directionally differentiable on* $C^{o}$ *such that, for any* $x \in C^{o}$, *we have*

$$\Pi'_C(x)(v) = v, \quad \text{for every } v \in X \text{ with } v \neq \theta.$$

*Proof*. For any given point $x \in C^{o}$, there is $\delta_x > 0$ such that

$$x + w \in C^{o}, \text{ for any } w \in X \text{ with } \|w\| < \delta_x.$$

This implies

$$\Pi_C(x) = x \text{ and } \Pi(x + w) = x + w, \text{ for any } w \in X \text{ with } \|w\| < \delta_x. \tag{4.2}$$

For any $v \in X$ with $v \neq \theta$, by (3.2), we calculate

$$\begin{aligned}\Pi'_C(x)(v) &= \lim_{t\downarrow 0}\frac{\Pi_C(x+tv)-\Pi_C(x)}{t}\\ &= \lim_{t\downarrow 0,\, t<\delta_x\backslash\|v\|}\frac{\Pi_C(x+tv)-x}{t}\\ &= \lim_{t\downarrow 0,\, t<\delta_x\backslash\|v\|}\frac{x+tv-x}{t}\\ &= \lim_{t\downarrow 0,\, t<\delta_x\backslash\|v\|}\frac{tv}{t}\\ &= v.\end{aligned}$$

□

## 4. Directional differentiability of the generalized metric projection Π onto closed balls

### 4.1. The generalized metric projection Π onto closed balls

For $r > 0$, the closed, open balls and the sphere in $X$ with radius $r$ and at center θ are respectively denoted by $B(r)$, $B^{o}(r)$ and $S(r)$. In particular, $B(X)$ denotes the closed unit ball and $S(X)$ denotes the unit sphere in $X$. As we reviewed in section 2, in general, $\Pi_C \neq P_C$. However, in particular, if $C = B(r)$ is a closed ball in $X$, then $\Pi_{B(r)}$ coincides with $P_{B(r)}$, which is proved by the following proposition and Lemma 3.5 in [13].

**Proposition 4.1** *Let* $r > 0$. *Then*

(a) *For any* $x \in B(r)$, *we have*

$$\Pi_{B(r)}(x) = x;$$

(b) *For any* $x \in X\backslash B(r)$, *we have*

$$\Pi_{B(r)}(x) = \frac{r}{\|x\|}x.$$

*Proof*. Part (a) is clear. We only prove part (b). For an arbitrary given $x \in X\backslash B(r)$, it is equivalent to $\|x\| > r$. By the semi-linear property of $J$ and by $1 - \frac{r}{\|x\|} > 0$, for any $z \in B(r)$ with $\|z\| \leq r$, we have

$$\begin{aligned}
&\langle J(x) - J(\tfrac{r}{\|x\|}x), \tfrac{r}{\|x\|}x - z\rangle \\
&= \left(1 - \tfrac{r}{\|x\|}\right) \langle J(x), \tfrac{r}{\|x\|}x - z\rangle \\
&= \left(1 - \tfrac{r}{\|x\|}\right) (r\|x\| - \langle J(x), z\rangle) \\
&\geq \left(1 - \tfrac{r}{\|x\|}\right) (r\|x\| - \|x\|\|z\|) \\
&\geq 0, \text{ for any } z \in B(r).
\end{aligned}$$

By the basic variational principle of $\Pi_{B(r)}$, this proves part (b). □

### 4.2 Directional differentiability of the generalized metric projection Π onto closed balls

Before we study the directional differentiability of the generalized metric projection onto closed balls, we need recall the following notations introduced in [13]. For any $x \in S(r)$, we define two subsets $x_r^{\uparrow}$ and $x_r^{\downarrow}$ in $X\backslash\{\theta\}$ as follows: for $v \in X$ with $v \neq \theta$, we say

(a) $v \in x_r^{\uparrow} \iff$ there is $\delta > 0$ such that $\|x + tv\| \geq r$, for all $t \in (0, \delta)$;
(b) $v \in x_r^{\downarrow} \iff$ there is $\delta > 0$ such that $\|x + tv\| < r$, for all $t \in (0, \delta)$.

In particular, for any given $x \in S(X)$ and for $v \in X$ with $v \neq \theta$, we write

(c) $x^{\uparrow} = x_1^{\uparrow}$: $v \in x^{\uparrow} \iff$ there is $\delta > 0$ such that $\|x + tv\| \geq 1$, for all $t \in (0, \delta)$;
(d) $x^{\downarrow} = x_1^{\downarrow}$: $v \in x^{\downarrow} \iff$ there is $\delta > 0$ such that $\|x + tv\| < 1$, for all $t \in (0, \delta)$.

The Lemma 5.1 in [13] shows that, for any given $x \in S(r)$, the two subsets $x_r^{\uparrow}$ and $x_r^{\downarrow}$ form a partition of $X\backslash\{\theta\}$. We list Lemma 5.1 in [13] below with slightly modified version.

**Lemma 5.1 in [13]** *Let $r > 0$. Then, for any $x \in S(r)$, we have*

$$x_r^{\uparrow} \cap x_r^{\downarrow} = \emptyset \quad \textit{and} \quad x_r^{\uparrow} \cup x_r^{\downarrow} = X\backslash\{\theta\}.$$

Theorem 5.2 in [13] proves that the metric projection on closed balls is directionally differentiable on the whole space $X$ and it gives the exact analytic representations (solutions) of the directional derivatives of the metric projection on closed balls. Since for any closed ball in $X$, by Proposition 4.1 in this paper and Lemma 3.5 in [13], we have that

$$\Pi_{B(r)} = P_{B(r)}, \text{ for any } r > 0. \tag{4.1}$$

It follows immediately that we must have

$$\Pi'_{B(r)}(x)(v) = P'_{B(r)}(x)(v), \text{ for every } x, v \in X \text{ with } v \neq \theta. \tag{4.2}$$

We show the details of (4.2) which are from Theorem 5.2 in [13].

**Theorem 4.2** *Let $B(r)$ be a closed ball in $X$ with $r > 0$. Then, $\Pi_{B(r)}$ is directionally differentiable on $X$ such that, for every $v \in X$ with $v \neq \theta$, we have*

(i) *For any $x \in (B(r))^{o}$,*

$$\Pi'_{B(r)}(x)(v) = v;$$

(ii) *For any* $x \in X\backslash B(r)$,

$$\Pi'_{B(r)}(x)(v) = \frac{r}{\|x\|^2}\left(\|x\|v - \Psi\left(\frac{x}{\|x\|}, \frac{v}{\|v\|}\right)\|v\|x\right).$$

*In particular,*

$$\Pi'_{B(r)}(x)(x) = \theta, \textit{ for any } x \in X\backslash B(r);$$

(iii) *For any* $x \in S(r)$,

(a)
$$\Pi'_{B(r)}(x)(v) = v - \frac{\|v\|}{r}\Psi\left(\frac{x}{r}, \frac{v}{\|v\|}\right)x, \quad \textit{if } v \in x_r^{\uparrow}.$$

*In particular,*

$$\Pi'_{B(r)}(x)(x) = \theta, \textit{ for any } x \in S(r);$$

(b)
$$\Pi'_{B(r)}(x)(v) = v, \quad \textit{if } v \in x_r^{\downarrow}.$$

*Proof*. By (4.2), the main parts of this theorem immediately follow from Proposition 4.1 in this paper and Theorem 5.2 in [13] (replace $c$ by $\theta$ in Theorem 5.2 in [13]). To complete the proof of (ii), by Lemma 2.2, we have

$$\Psi\left(\frac{x}{\|x\|}, \frac{x}{\|x\|}\right) = \left\|\frac{x}{\|x\|}\right\| = 1, \text{ for any } x \in X\backslash B(r).$$

Notice that

$$x \in x_r^{\uparrow}, \text{ for any } x \in S(r).$$ □

## 5. Directional differentiability of the generalized metric projection Π onto subspaces and closed and convex cones

### 5.1. The generalized metric projection Π onto closed and convex cones

In this section, we study the directional differentiability of the generalized metric projection on subspaces and closed and convex cones in a uniformly convex and uniformly smooth Banach space $X$, which has dual space $X^*$. We also give the analytic representations of the directional derivatives of the generalized metric projection. In the first subsection, we study the solutions of the generalized metric projection on subspaces and closed and convex cones. These solutions will help us to find the directional derivatives of the generalized metric projection onto such special closed and convex subsets of $X$.

**Lemma 5.1**. *Let $K$ be a closed and convex cone with vertex at* θ *in $X$. For $x \in X$, we have*

$$\Pi_K(\lambda x) = \lambda\Pi_K(x), \textit{ for any } \lambda > 0.$$

*Proof*. By the variational principle of $\Pi_K$, for $x \in X$, we have

$$\langle J(x) - J(\Pi_K(x)), \Pi_K(x) - z\rangle \geq 0, \text{ for any } z \in K.$$

By the semi-linearity of $J$, it follows that

$$\begin{aligned}&\langle J(\lambda x) - J(\lambda\Pi_K(x)), \lambda\Pi_K(x) - z\rangle\\ &= \lambda^2\langle J(x) - J(\Pi_K(x)), \Pi_K(x) - \tfrac{1}{\lambda}z\rangle\\ &\geq 0, \text{ for any } z \in K.\end{aligned}$$

By the variational principle of $\Pi_C$ again, this implies $\Pi_K(\lambda x) = \lambda\Pi_K(x)$. □

**Lemma 5.2**. *Let $K$ be a closed and convex cone in $X$ with vertex at $c \in X$. For $x \in X$, we have*

$$\langle J(x) - J(\Pi_K(x)), \Pi_K(x) - c\rangle = 0.$$

*Proof*. It is clear for $\Pi_K(x) = c$. Next, we suppose that $\Pi_K(x) \neq c$. For $x \in X$, by $c \in K$ and by the variational principle of $\Pi_C$, we have

$$\langle J(x) - J(\Pi_K(x)), \Pi_K(x) - c\rangle \geq 0.$$

Take $K \ni z = c + 2(\Pi_K(x) - c) = 2\Pi_K(x) - c$, by the variational principle of $\Pi_C$ again, we have

$$\begin{aligned} &\langle J(x) - J(\Pi_K(x)), -\Pi_K(x) + c\rangle \\ &= \langle J(x) - J(\Pi_K(x)), \Pi_K(x) - z\rangle \\ &\geq 0. \end{aligned}$$

This proves this lemma. □

For closed subspaces of $X$, the basic variational principle of the generalized metric projection become the following orthogonal version.

**Lemma 5.3**. *Let $Y$ be a closed subspace of $X$. For $x \in X$ and $y \in Y$, we have*

$$y = \Pi_Y(x) \quad \Leftrightarrow \quad \langle J(x) - J(y), z\rangle = 0, \text{ for all } z \in Y.$$

*Proof*. Let $x \in X$ and $y \in Y$. For any $z \in Y$, since $y - z \in Y$, by the variational principle of $\Pi_Y$, we have that $y = \Pi_Y(x)$ if and only if

$$\begin{aligned} &\langle J(x) - J(y), z\rangle \\ &= \langle J(x) - J(y), y - (y - z)\rangle \\ &\geq 0. \end{aligned}$$

By $y + z \in Y$, we have that $y = \Pi_Y(x)$ if and only if

$$\begin{aligned} &\langle J(x) - J(y), -z\rangle \\ &= \langle J(x) - J(y), y - (y + z)\rangle \\ &\geq 0. \end{aligned}$$

This proves this lemma. □

Let $C$ be a nonempty closed and convex subset of $X$. By Definition 3.8 in [12], the orthogonal cone (in $X$) of $C$ is

$$C^{\perp} = \{x \in X: \langle Jx, z\rangle = 0, \text{ for all } z \in C\}.$$

The dual orthogonal cone in $X^*$ of $C$ is

$$C^{\llcorner} = \{\varphi \in X^*: \langle \varphi, z\rangle = 0, \text{ for all } z \in C\}.$$

By Corollary 2.8 in [12], we have the properties of the orthogonal cones of $C$.

**Lemma 5.4**. *Let $C$ be a nonempty closed and convex subset of $X$. We have*

(a) $C^{\perp}$ *is a closed cone in $X$ with vertex at* $\theta$, *which is not convex, in general*;
(b) $C^{\llcorner}$ *is a closed subspace in* $X^*$.

*Proof*. Part (a) is from Corollary 2.8 in [12]. The proof of (b) is straight forward and is omitted here. □

In particular, by the dual orthogonal cones in $X^*$ of closed subspaces in $X$, we have the following results.

**Corollary 5.5**. *Let Y be a closed subspace of X. For any y ∈ Y, we have*

$$\Pi_Y^{-1}(y) = \{x \in X: J(x) - J(y) \in Y^{\llcorner}\} = \{x \in X: J(x) \in J(y) + Y^{\llcorner}\}.$$

*Proof*. For any $x \in X$ and $y \in Y$, by Lemma 5.3 and by the definition of $Y^{\llcorner}$, we have

$$\Pi_Y(x) = y \quad \Leftrightarrow \quad J(x) - J(y) \in Y^{\llcorner}.$$

It is equivalent to

$$J(x) \in J(y) + Y^{\llcorner}.$$ □

### 5.2. Directional differentiability of the generalized metric projection Π onto subspaces and closed and convex cones

By Lemma 3.2, we have the following results immediately.

**Corollary 5.6** *Let $c \in X$ and let K be a closed and convex cone in X with vertex at c. Then, for any $y \in K$ with $y \neq c$, $\Pi_K$ is directionally differentiable at c along direction $y - c$ such that*

$$\Pi_K'(c; y - c) = y - c.$$

**Lemma 5.7**. *Let K be a closed and convex cone in X with vertex at θ. We have*

(a) $\Pi_K'(x; x) = \Pi_K(x)$, *for any* $x \in X\backslash\{\theta\}$;
(b) $\Pi_K'(x; x) = x$, *for any* $x \in K\backslash\{\theta\}$;
(c) $\Pi_K'(\theta; x) = x$, *for any* $x \in K\backslash\{\theta\}$.

*Proof*. It is clear that part (b) follows from part (a); and part (c) follows from Lemma 5.6 immediately. So, we only need to show part (a). For any $x \in X\backslash\{\theta\}$, by Lemma 5.1, we calculate

$$\begin{aligned}\Pi_K'(x; x) &= \lim_{t\downarrow 0}\frac{\Pi_K(x+tx)-\Pi_K(x)}{t}\\ &= \lim_{t\downarrow 0}\frac{(1+t)\Pi_K(x)-\Pi_K(x)}{t}\\ &= \Pi_K(x).\end{aligned}$$

This proves (a). □

**Lemma 5.8**. *Let Y be a closed subspace in X. For any y ∈ Y, we have*

$$\Pi_Y'(y; v) = v, \text{ for any } v \in Y\backslash\{\theta\}.$$

Proof. Taking arbitrarily $y \in Y$ and $v \in Y\backslash\{\theta\}$, by $Y$ being a subspace of $X$, we have

$$\begin{aligned}\Pi_Y'(y; v) &= \lim_{t\downarrow 0}\frac{\Pi_Y(y+tv)-\Pi_Y(y)}{t}\\ &= \lim_{t\downarrow 0}\frac{y+tv-y}{t}\\ &= v.\end{aligned}$$ □

## 6. Real $l_p$ space

### 6.1. The generalized metric projection Π onto closed balls in real $l_p$ spaces

Through this section, let $p$ be a positive number with $1 < p < \infty$ and let $(l_p, \|\cdot\|_{l_p})$ be the Banach space of real sequences satisfying

$$\|x\|_{l_p} = \left(\sum_{n=1}^{\infty}|x_n|^p\right)^{\frac{1}{p}} < \infty, \text{ for any } x = (x_1, x_2, \ldots) \in l_p.$$

$(l_p, \|\cdot\|_{l_p})$ is a special example of uniformly convex and uniformly smooth Banach spaces with dual space $(l_q, \|\cdot\|_{l_q})$, where $q$ is a positive number with $1 < q < \infty$ satisfying $\frac{1}{p} + \frac{1}{q} = 1$. Recall that the normalized duality mappings $J$ in $l_p$ holds the following analytic representations. For any $x = (x_1, x_2, \ldots) \in l_p$ with $x \neq \theta$, we have $Jx = ((Jx)_1, (Jx)_2, \ldots) \in (l_p)^* = l_q$, such that

$$(Jx)_n = \frac{|x_n|^{p-1}\text{sign}(x_n)}{\|x\|_{l_p}^{p-2}} = \frac{|x_n|^{p-2}x_n}{\|x\|_{l_p}^{p-2}}, \text{ for } n = 1, 2, \ldots. \tag{6.1}$$

Let $\mathbb{N}$ denote the set of all positive integers. In this section, let $M$ be a nonempty subset of $\mathbb{N}$. Let $l_p^M$ denote the following subset of $l_p$:

$$l_p^M = \{x = (x_1, x_2, \ldots) \in l_p \colon x_n = 0, \text{ for } n \notin M\}.$$

Then, $(l_p^M, \|\cdot\|_{l_p})$ is a closed subspace in $l_p$ such that, for any $x = (x_1, x_2, \ldots) \in l_p^M$, we have

$$\|x\|_{l_p} = \left(\sum_{m\in M}|x_m|^p\right)^{\frac{1}{p}}.$$

If $\mathbb{N}\backslash M \neq \emptyset$, then, $l_p^M$ is a proper closed subspace of $l_p$. In the case $\mathbb{N}\backslash M \neq \emptyset$, $l_p^{\mathbb{N}\backslash M}$ is similarly defined to be a proper closed subspace of $l_p$. For any $x = (x_1, x_2, \ldots) \in l_p$, we define $x^M \in l_p^M$ by

$$(x^M)_n = \begin{cases} x_n, \text{for } n \in M, \\ 0, \text{ for } n \notin M. \end{cases} \tag{6.2}$$

It is well-known that the normalized duality mapping $J$ is not linear on $l_p$. That is, for $x, y \in l_p$, $J(x + y) \neq J(x) + J(y)$, in general. However, by (6.1), when $J$ is restricted on $l_p^M$, it has the following useful properties.

**Lemma 6.1**. *The normalized duality mapping J holds the following properties*:

(a) *J maps points in $l_p^M$ to points in $l_p^M$. That is,*

$$J(x) \in l_p^M, \text{ for any } x \in l_p^M.$$

(b) *For any $x \in l_p\backslash\{\theta\}$, we have*

$$J(x) = \frac{\|x^M\|_{l_p}^{p-2}}{\|x\|_{l_p}^{p-2}}J(x^M) + \frac{\|x^{\mathbb{N}\backslash M}\|_{l_p}^{p-2}}{\|x\|_{l_p}^{p-2}}J(x^{\mathbb{N}\backslash M}), \tag{6.3}$$

*and*

$$\|x\|_{l_p}^p = \|x^M\|_{l_p}^p + \|x^{\mathbb{N}\backslash M}\|_{l_p}^p.$$

*Proof*. Proof of (a). Let $x = (x_1, x_2, \ldots) \in l_p^M$. Then, for any $n \notin M$, we have $x_n = 0$. It follows that

$$(J(x))_n = \frac{|x_n|^{p-2}x_n}{\|x\|_{l_p}^{p-2}} = 0, \text{ for any } n \notin M.$$

This implies

$$J(x) \in l_p^M, \text{ for any } x \in l_p^M.$$

Proof of (b). For any $x = (x_1, x_2, \ldots) \in l_p$, by (a), it is clear that this lemma holds if it satisfies any one of the following equations: $x = \theta$, $x^M = \theta$, or $x^{\mathbb{N}\backslash M} = \theta$. So, we suppose $x \neq \theta$, $x^M \neq \theta$, and $x^{\mathbb{N}\backslash M} \neq \theta$. By (6.1), for $n \in M$, which satisfies $(x^{\mathbb{N}\backslash M})_n = 0$, we have

$$\begin{aligned}
(J(x))_n &= \frac{|x_n|^{p-2}x_n}{\|x\|_{l_p}^{p-2}} \\
&= \frac{\|x^M\|_{l_p}^{p-2}}{\|x\|_{l_p}^{p-2}} \frac{|(x^M)_n|^{p-2}(x^M)_n}{\|x^M\|_{l_p}^{p-2}} \\
&= \frac{\|x^M\|_{l_p}^{p-2}}{\|x\|_{l_p}^{p-2}} \frac{|(x^M)_n|^{p-2}(x^M)_n}{\|x^M\|_{l_p}^{p-2}} + \frac{\|x^{\mathbb{N}\backslash M}\|_{l_p}^{p-2}}{\|x\|_{l_p}^{p-2}} \frac{|(x^{\mathbb{N}\backslash M})_n|^{p-2}(x^{\mathbb{N}\backslash M})_n}{\|x^{\mathbb{N}\backslash M}\|_{l_p}^{p-2}}. \\
&= \frac{\|x^M\|_{l_p}^{p-2}}{\|x\|_{l_p}^{p-2}}(J(x^M))_n + \frac{\|x^{\mathbb{N}\backslash M}\|_{l_p}^{p-2}}{\|x\|_{l_p}^{p-2}}(J(x^{\mathbb{N}\backslash M}))_n.
\end{aligned}$$

We similarly prove the above equation for any $n \in \mathbb{N}\backslash M$. □

For any $r > 0$, the closed, open balls and the sphere in $l_p^M$ with radius $r$ and with center at the origin are respectively denoted by

$$B_M(r) = \left\{x \in l_p^M : \|x\|_{l_p} = \left(\sum_{m\in M}|x_m|^p\right)^{\frac{1}{p}} \leq r\right\},$$

$$B_M^o(r) = \left\{x \in l_p^M : \|x\|_{l_p} = \left(\sum_{m\in M}|x_m|^p\right)^{\frac{1}{p}} < r\right\},$$

$$S_M(r) = \left\{x \in l_p^M : \|x\|_{l_p} = \left(\sum_{m\in M}|x_m|^p\right)^{\frac{1}{p}} = r\right\}.$$

$B_M(r)$ is a nonempty closed, bounded and convex subset in $l_p$. However, if $\mathbb{N}\backslash M \neq \emptyset$, then $B_M(r)$ is not a (closed) ball in $l_p$. We define the following cylinders in $l_p$:

$$C_M(r) = \{x \in l_p : x^M \in B_M(r)\},$$

$$C_M^o(r) = \{x \in l_p : x^M \in B_M^o(r)\}.$$

$C_M(r)$ and $C_M^o(r)$ are convex cylinders in $l_p$ with bases $B_M(r)$ and $B_M^o(r)$ in $l_p^M$, respectively. It is clear to see that if $\mathbb{N}\backslash M \neq \emptyset$, then $C_M(r)$ is a closed, unbounded and convex subset in $l_p$ and $C_A^o(r)$ is an open, unbounded and convex subset in $l_p$. In particular, if $\mathbb{N}\backslash M = \emptyset$, then $l_p^M$ coincides with $l_p$. In this case, $B_M(r)$, $B_M^o(r)$, $S_M(r)$, $C_M(r)$ and $C_M^o(r)$ are respectively denoted by $B(r)$, $B^o(r)$, $S(r)$, $C(r)$ and $C^o(r)$, which satisfy

$$B(r) = C(r) \quad \text{and} \quad B^o(r) = C^o(r).$$

Next, we calculate the values of the generalized metric projection $\Pi_{B_M(r)}$ from $l_p$ onto $B_M(r)$.

**Theorem 6.2**. *For any $r > 0$ and $x \in l_p$, we have*

(a) $\Pi_{B_M(r)}(x) = x$, *for* $x \in B_M(r)$;

(b) $\Pi_{B_M(r)}(x) = \frac{r}{\|x\|_{l_p}} x$, *for* $x \in l_p^M \backslash B_M(r)$;

(c) $\Pi_{B_M(r)}(x) = \frac{\|x^M\|_{l_p}^{p-2}}{\|x\|_{l_p}^{p-2}} x^M$, *for* $x \in C_M(r) \backslash l_p^M$;

(d) $\Pi_{B_M(r)}(x) = \begin{cases} \frac{r}{\|x^M\|_{l_p}} x^M, & if \ \frac{\|x^M\|_{l_p}^{p-2}}{\|x\|_{l_p}^{p-2}} > \frac{r}{\|x^M\|_{l_p}}, \\ \frac{\|x^M\|_{l_p}^{p-2}}{\|x\|_{l_p}^{p-2}} x^M, & if \ \frac{\|x^M\|_{l_p}^{p-2}}{\|x\|_{l_p}^{p-2}} \leq \frac{r}{\|x^M\|_{l_p}}, \end{cases}$ *for* $x \in l_p \backslash (C_M(r) \cup l_p^M)$.

*Proof*. The proof of this theorem is similar to the proof of Theorem 6.1 in [13]. Part (a) is clear. We prove part (b). For any given $x \in l_p^M \backslash B_M(r)$, $x$ must satisfies $\|x\|_{l_p} > r$. Then, for any $y \in B_M(r)$, we have

$$\begin{aligned} &\langle J(x) - J\left(\frac{r}{\|x\|_{l_p}} x\right), \ \frac{r}{\|x\|_{l_p}} x - y\rangle \\ &= \left(1 - \frac{r}{\|x\|_{l_p}}\right)\langle J(x), \ \frac{r}{\|x\|_{l_p}} x - y\rangle \\ &= \left(1 - \frac{r}{\|x\|_{l_p}}\right)\left(r\|x\|_{l_p} - \langle J(x), y\rangle\right) \\ &\geq \left(1 - \frac{r}{\|x\|_{l_p}}\right)\left(r\|x\|_{l_p} - \|J(x)\|_{l_q}\|y\|_{l_p}\right) \\ &= \left(1 - \frac{r}{\|x\|_{l_p}}\right)\left(r\|x\|_{l_p} - \|x\|_{l_p}\|y\|_{l_p}\right) \\ &= \left(\|x\|_{l_p} - r\right)\left(r - \|y\|_{l_p}\right) \\ &\geq 0, \text{ for all } y \in B_M(r). \end{aligned}$$

By $x \in l_p^M$ and $\|x\|_{l_p} > r$, we have $\frac{r}{\|x\|_{l_p}} x \in S_M(r)$. By the basic variational principle of $\Pi_{B_M(r)}$, the above inequality implies

$$\Pi_{B_M(r)}(x) = \frac{r}{\|x\|_{l_p}} x, \quad \text{for any } x \in l_p^M \backslash B_M(r).$$

Proof of (c). For $x \in C_M(r) \backslash l_p^M$, we have $x^M \in B_M(r)$ and $x \notin l_p^M$. It follows that

$$\|x^M\|_{l_p} \leq r \quad \text{and} \quad \|x^M\|_{l_p} < \|x\|_{l_p}.$$

Then, we have

$$\frac{\|x^M\|_{l_p}^{p-2}}{\|x\|_{l_p}^{p-2}} x^M \in B_M(r) \quad \text{and} \quad \frac{\|x^M\|_{l_p}^{p-2}}{\|x\|_{l_p}^{p-2}} x^M - y \in l_p^M, \text{ for any } y \in B_M(r) \subseteq l_p^M. \tag{6.4}$$

By (6.2), (6.3) and (6.4), we obtain

$$\begin{aligned} &\langle J(x) - J\left(\frac{\|x^M\|_{l_p}^{p-2}}{\|x\|_{l_p}^{p-2}} x^M\right), \ \frac{\|x^M\|_{l_p}^{p-2}}{\|x\|_{l_p}^{p-2}} x^M - y\rangle \\ &= \langle \frac{\|x^{\mathbb{N}\backslash M}\|_{l_p}^{p-2}}{\|x\|_{l_p}^{p-2}} J(x^{\mathbb{N}\backslash M}) + \frac{\|x^M\|_{l_p}^{p-2}}{\|x\|_{l_p}^{p-2}} J(x^M) - J\left(\frac{\|x^M\|_{l_p}^{p-2}}{\|x\|_{l_p}^{p-2}} x^M\right), \ \frac{\|x^M\|_{l_p}^{p-2}}{\|x\|_{l_p}^{p-2}} x^M - y\rangle \end{aligned}$$

$$= \langle \frac{\|x^{\mathbb{N}\setminus M}\|_{l_p}^{p-2}}{\|x\|_{l_p}^{p-2}} J(x^{\mathbb{N}\setminus M}) + \frac{\|x^M\|_{l_p}^{p-2}}{\|x\|_{l_p}^{p-2}} J(x^M) - \frac{\|x^M\|_{l_p}^{p-2}}{\|x\|_{l_p}^{p-2}} J(x^M),\ \frac{\|x^M\|_{l_p}^{p-2}}{\|x\|_{l_p}^{p-2}} x^M - y\rangle$$
$$= \langle \frac{\|x^{\mathbb{N}\setminus M}\|_{l_p}^{p-2}}{\|x\|_{l_p}^{p-2}} J(x^{\mathbb{N}\setminus M}),\ \frac{\|x^M\|_{l_p}^{p-2}}{\|x\|_{l_p}^{p-2}} x^M - y\rangle$$
$$= 0, \text{ for all } y \in B_M(r).$$

By (6.4) again and by the basic variational principle of $\Pi_{B_M(r)}$, this proves part (c). Next, we prove part (d). For any $x \in l_p \setminus (C_M(r) \cup l_p^M)$, it follows that

$$r < \|x^M\|_{l_p} < \|x\|_{l_p}. \tag{6.5}$$

Case 1. $\frac{\|x^M\|_{l_p}^{p-2}}{\|x\|_{l_p}^{p-2}} > \frac{r}{\|x^M\|_{l_p}}$. In this case, we consider the point $\frac{r}{\|x^M\|_{l_p}} x^M$. By (6.5), we have

$$\frac{r}{\|x^M\|_{l_p}} x^M \in S_M(r) \subseteq B_M(r), \text{ for any } x \in l_p \setminus (C_M(r) \cup l_p^M) \text{ with } \frac{\|x^M\|_{l_p}^{p-2}}{\|x\|_{l_p}^{p-2}} > \frac{r}{\|x^M\|_{l_p}}. \tag{6.6}$$

We also have

$$\frac{r}{\|x^M\|_{l_p}} x^M - y \in l_p^M, \text{ for any } y \in B_M(r) \subseteq l_p^M. \tag{6.7}$$

Then, for any $y \in B_M(r) \subseteq l_p^M$, by (6.3) in Lemma 6.1 and (6.7), in Case 1, we obtain

$$\langle J(x) - J(\frac{r}{\|x^M\|_{l_p}} x^M),\ \frac{r}{\|x^M\|_{l_p}} x^M - y\rangle$$
$$= \langle \frac{\|x^{\mathbb{N}\setminus M}\|_{l_p}^{p-2}}{\|x\|_{l_p}^{p-2}} J(x^{\mathbb{N}\setminus M}) + \frac{\|x^M\|_{l_p}^{p-2}}{\|x\|_{l_p}^{p-2}} J(x^M) - \frac{r}{\|x^M\|_{l_p}} J(x^M),\ \frac{r}{\|x^M\|_{l_p}} x^M - y\rangle$$
$$= \langle \frac{\|x^{\mathbb{N}\setminus M}\|_{l_p}^{p-2}}{\|x\|_{l_p}^{p-2}} J(x^{\mathbb{N}\setminus M}), \frac{r}{\|x^M\|_{l_p}} x^M - y\rangle + \langle \left(\frac{\|x^M\|_{l_p}^{p-2}}{\|x\|_{l_p}^{p-2}} - \frac{r}{\|x^M\|_{l_p}}\right) J(x^M),\ \frac{r}{\|x^M\|_{l_p}} x^M - y\rangle$$
$$= \left(\frac{\|x^M\|_{l_p}^{p-2}}{\|x\|_{l_p}^{p-2}} - \frac{r}{\|x^M\|_{l_p}}\right) \langle J(x^M), \frac{r}{\|x^M\|_{l_p}} x^M - y\rangle$$
$$= \left(\frac{\|x^M\|_{l_p}^{p-2}}{\|x\|_{l_p}^{p-2}} - \frac{r}{\|x^M\|_{l_p}}\right) \left(\frac{r}{\|x^M\|_{l_p}} \|x^M\|_{l_p}^2 - \langle J(x^M), y\rangle\right)$$
$$\geq \left(\frac{\|x^M\|_{l_p}^{p-2}}{\|x\|_{l_p}^{p-2}} - \frac{r}{\|x^M\|_{l_p}}\right) \left(r\|x^M\|_{l_p} - \|x^M\|_{l_p}\|y\|_{l_p}\right)$$
$$= \left(\frac{\|x^M\|_{l_p}^{p-2}}{\|x\|_{l_p}^{p-2}} - \frac{r}{\|x^M\|_{l_p}}\right) \|x^M\|_{l_p} \left(r - \|y\|_{l_p}\right)$$
$$\geq 0, \text{ for all } y \in B_M(r).$$

By (6.6) and the basic variational principle of $\Pi_{B_M(r)}$, this implies that, for $x \in l_p \setminus (C_M(r) \cup l_p^M)$,

$$\Pi_{B_M(r)}(x) = \frac{r}{\|x^M\|_{l_p}} x^M, \quad \text{if } \frac{\|x^M\|_{l_p}^{p-2}}{\|x\|_{l_p}^{p-2}} > \frac{r}{\|x^M\|_{l_p}}.$$

Case 2. $\frac{\|x^M\|_{l_p}^{p-2}}{\|x\|_{l_p}^{p-2}} \leq \frac{r}{\|x^M\|_{l_p}}$. In this case, we consider the point $\frac{\|x^M\|_{l_p}^{p-2}}{\|x\|_{l_p}^{p-2}} x^M$. It is clear to see that

$$\left\| \frac{\|x^M\|_{l_p}^{p-2}}{\|x\|_{l_p}^{p-2}} x^M \right\|_{l_p} \leq \left\| \frac{r}{\|x^M\|_{l_p}} x^M \right\|_{l_p} = r.$$

This implies that

$$\frac{\|x^M\|_{l_p}^{p-2}}{\|x\|_{l_p}^{p-2}} x^M \in B_M(r), \text{ for any } x \in l_p \backslash (C_M(r) \cup l_p^M) \text{ with } \frac{\|x^M\|_{l_p}^{p-2}}{\|x\|_{l_p}^{p-2}} \leq \frac{r}{\|x^M\|_{l_p}}. \tag{6.8}$$

Then, for any $y \in B_M(r) \subseteq l_p^M$, by (6.3) in Lemma 6.1 and (6.8) in Case 2, we obtain

$$\begin{aligned}
&\langle J(x) - J(\frac{\|x^M\|_{l_p}^{p-2}}{\|x\|_{l_p}^{p-2}} x^M), \frac{\|x^M\|_{l_p}^{p-2}}{\|x\|_{l_p}^{p-2}} x^M - y \rangle \\
&= \langle \frac{\|x^{\mathbb{N}\backslash M}\|_{l_p}^{p-2}}{\|x\|_{l_p}^{p-2}} J(x^{\mathbb{N}\backslash M}) + \frac{\|x^M\|_{l_p}^{p-2}}{\|x\|_{l_p}^{p-2}} J(x^M) - \frac{\|x^M\|_{l_p}^{p-2}}{\|x\|_{l_p}^{p-2}} J(x^M), \frac{\|x^M\|_{l_p}^{p-2}}{\|x\|_{l_p}^{p-2}} x^M - y \rangle \\
&= \langle \frac{\|x^{\mathbb{N}\backslash M}\|_{l_p}^{p-2}}{\|x\|_{l_p}^{p-2}} J(x^{\mathbb{N}\backslash M}), \frac{\|x^M\|_{l_p}^{p-2}}{\|x\|_{l_p}^{p-2}} x^M - y \rangle \\
&= 0, \text{ for all } y \in B_M(r).
\end{aligned}$$

By (6.8) and the basic variational principle of $\Pi_{B_M(r)}$, this implies that, for $x \in l_p \backslash (C_M(r) \cup l_p^M)$,

$$\Pi_{B_M(r)}(x) = \frac{\|x^M\|_{l_p}^{p-2}}{\|x\|_{l_p}^{p-2}} x^M, \quad \text{if } \frac{\|x^M\|_{l_p}^{p-2}}{\|x\|_{l_p}^{p-2}} \leq \frac{r}{\|x^M\|_{l_p}}. \qquad \square$$

Notice that if $\mathbb{N}\backslash M = \emptyset$, then $l_p^M$ coincides with $l_p$. In this case, we have

$$C_M(r) = B_M(r) = B(r) \quad \text{and} \quad x^M = x, \text{ for any } x \in l_p.$$

Then, by Theorem 6.2, we obtain the following results.

**Corollary 6.3**. *For any $r > 0$ and $x \in l_p$, we have*

(a) $\Pi_{B(r)}(x) = x$, *for* $x \in B(r)$;
(b) $\Pi_{B(r)}(x) = \frac{r}{\|x\|_{l_p}} x$, *for* $x \in l_p \backslash B(r)$.

For $p = 2$, in the Hilbert space $l_2$, the generalized matric projection $\Pi$ coincides with the metric projection $P$, that is, $\Pi_{B_M(r)} = P_{B_M(r)}$. From Theorem 6.2, we obtain

**Corollary 6.4**. *For any $r > 0$ and $x \in l_2$, we have*

(a) $\Pi_{B_M(r)}(x) = P_{B_M(r)}(x) = x$, *for* $x \in B_M(r)$;

(b) $\Pi_{B_M(r)}(x) = P_{B_M(r)}(x) = \frac{r}{\|x\|_{l_2}} x,$ *for* $x \in l_2^M \backslash B_M(r)$.

(c) $\Pi_{B_M(r)}(x) = P_{B_M(r)}(x) = x^M,$ *for* $x \in C_M(r) \backslash l_2^M$;

(d) $$\Pi_{B_M(r)}(x) = P_{B_M(r)}(x) = \begin{cases} \frac{r}{\|x^M\|_{l_2}} x^M, & if \ \frac{r}{\|x^M\|_{l_2}} < 1, \\ x^M, & if \ \frac{r}{\|x^M\|_{l_2}} \geq 1, \end{cases} \quad for\ x \in l_2 \backslash (C_M(r) \cup l_2^M).$$

*In particular, if* $\mathbb{N} \backslash M = \emptyset$, *then,*

(e) $\Pi_{B(r)}(x) = P_{B(r)}(x) = x,$ *for* $x \in B(r)$;

(f) $\Pi_{B(r)}(x) = P_{B(r)}(x) = \frac{r}{\|x\|_{l_2}} x,$ *for* $x \in l_2^M \backslash B(r)$.

### 6.2. The directional differentiability of the generalized metric projection Π onto closed balls in $l_p$

By the solutions of the generalized metric projection onto $B(r)$ in Theorem 6.2, in this subsection, we study the directional differentiability of $\Pi_{B(r)}$. We will see that the solutions of $\Pi'_{B_M(r)}$ are much more complicated than the solutions of $P'_{B_M(r)}$.

**Theorem 6.5**. *For any* $r > 0$ *and for any* $x, h \in l_p$ *with* $h \neq \theta$, *we have*

(a) *For* $x \in B_M^o(r)$,

(i) $\Pi'_{B_M(r)}(x)(h) = h,$ *for* $h \in l_p^M$;

(ii) *If p is a positive integer greater than* 1, *then,*

$$\Pi'_{B_M(r)}(x)(h) = \frac{\|x^M\|_{l_p}^{p-2}}{\|x\|_{l_p}^{p-2}} h^M + \frac{(p-2)}{\|x\|_{l_p}} \big(\Psi(x, h^M) - \Psi(x, h)\big)x, \ for\ h \notin l_p^M.$$

*In particular, for* $p = 2$,

$$\Pi'_{B_M(r)}(x)(h) = h^M, \text{ for } h \in l_2 \backslash l_2^M;$$

(b) *If* $x \in l_p^M \backslash B_M(r)$, *then*

$$\Pi'_{B_M(r)}(x)(h) = \frac{r}{\|x\|_{l_p}} h^M - \frac{r}{\|x\|_{l_p}^2} \Psi(x, h^M)x;$$

*In particular, we have*

(i) $\Pi'_{B_M(r)}(x)(h) = \frac{r}{\|x\|_{l_p}} h - \frac{r}{\|x\|_{l_p}^2} \Psi(x, h)x,$ *for* $h \in l_p^M$;

(ii) $\Pi'_{B_M(r)}(x)(x) = \theta$, *for any* $x \in l_p^M \backslash B_M(r)$;

(c) *If* $x \in C_M^o(r) \backslash l_p^M$, *then*

$$\Pi'_{B_M(r)}(x)(h) = h^M;$$

(d) *For* $x \in l_p \backslash (C_M(r) \cup l_p^M)$,

(i) *if x satisfies*

$$\frac{\|x^M\|_{l_p}^{p-2}}{\|x\|_{l_p}^{p-2}} > \frac{r}{\|x^M\|_{l_p}}, \tag{6.9}$$

*then*

$$\Pi'_{B_M(r)}(x)(h) = \frac{r}{\|x^M\|_{l_p}} h^M - \frac{r}{\|x^M\|_{l_p}^2} \Psi(x^M, h^M) x^M.$$

*In particular, for any* $x \in l_p \backslash (C_M(r) \cup l_p^M)$ *satisfying* (6.9), *then*

$$\Pi'_{B_M(r)}(x)(x) = \theta;$$

(ii) *if p is a positive integer greater than* 1 *and x satisfies*

$$\frac{\|x^M\|_{l_p}^{p-2}}{\|x\|_{l_p}^{p-2}} < \frac{r}{\|x^M\|_{l_p}}, \tag{6.10}$$

*then*

$$\Pi'_{B_M(r)}(x)(h) = \frac{\|x^M\|_{l_p}^{p-2}}{\|x\|_{l_p}^{p-2}} h^M + \frac{(p-2)\|x^M\|_{l_p}^{p-2}}{\|x\|_{l_p}^{p-2}} \left( \frac{\Psi(x^M, h^M)}{\|x^M\|_{l_p}} - \frac{\Psi(x,h)}{\|x\|_{l_p}} \right) x^M.$$

*In particular, for any* $x \in l_p \backslash (C_M(r) \cup l_p^M)$ *satisfying* (6.10), *then*

$$\Pi'_{B_M(r)}(x)(x) = x^M.$$

*Proof*. Proof of (a). Suppose $x \in B_M^o(r)$. If $h \in l_p^M$, then, there is $\delta > 0$, such that $x + th \in B_M^o(r)$, for all $t \in (0, \delta)$. By part (i) in (a) in Theorem 6.2, we have

$$\begin{aligned}
&\Pi'_{B_M(r)}(x)(h) \\
&= \lim_{t\downarrow 0} \frac{\Pi_{B_M(r)}(x+th) - \Pi_{B_M(r)}(x)}{t} \\
&= \lim_{t\downarrow 0, t<\delta} \frac{x+th-x}{t} \\
&= \lim_{t\downarrow 0, t<\delta} \frac{th}{t} \\
&= h.
\end{aligned}$$

Proof of part (ii) in (a) for $p = 2$. Let $x \in B_M^o(r)$. For $h \in l_2 \backslash l_2^M$, there is $\lambda > 0$, such that $x + th \in C_M(r) \backslash l_2^M$, for all $t \in (0, \lambda)$. By parts (a, c) in Corollary 6.4 and by (6.3), we have

$$\begin{aligned}
&\Pi'_{B_M(r)}(x)(h) \\
&= \lim_{t\downarrow 0} \frac{\Pi_{B_M(r)}(x+th) - \Pi_{B_M(r)}(x)}{t} \\
&= \lim_{t\downarrow 0, t<\lambda} \frac{(x+th)^M - x}{t} \\
&= \lim_{t\downarrow 0, t<\lambda} \frac{x+th^M - x}{t}
\end{aligned}$$

$$= \lim_{t\downarrow 0, t<\lambda} \frac{th^M}{t}$$
$$= h^M.$$

Proof of part (ii) in (a) for $p > 2$. Let $x \in B_M^o(r)$. For $h \in l_p \backslash l_p^M$, there is $\lambda > 0$, such that $x + th \in C_M(r) \backslash l_p^M$, for all $t \in (0, \lambda)$. By parts (a, c) in Theorem 6.2 and by (6.3), we have

$$\Pi'_{B_M(r)}(x)(h)$$
$$= \lim_{t\downarrow 0} \frac{\Pi_{B_M(r)}(x+th) - \Pi_{B_M(r)}(x)}{t}$$
$$= \lim_{t\downarrow 0, t<\lambda} \frac{\frac{\|(x+th)^M\|_{l_p}^{p-2}}{\|x+th\|_{l_p}^{p-2}}(x+th)^M - x}{t}$$
$$= \lim_{t\downarrow 0, t<\lambda} \frac{\frac{\|(x+th)^M\|_{l_p}^{p-2}}{\|x+th\|_{l_p}^{p-2}}(x+th^M) - x}{t}$$
$$= \lim_{t\downarrow 0, t<\lambda} \frac{\frac{\|(x+th)^M\|_{l_p}^{p-2}}{\|x+th\|_{l_p}^{p-2}} th^M + \left(\frac{\|(x+th)^M\|_{l_p}^{p-2}}{\|x+th\|_{l_p}^{p-2}} - 1\right)x}{t}$$
$$= \frac{\|x^M\|_{l_p}^{p-2}}{\|x\|_{l_p}^{p-2}} h^M + \lim_{t\downarrow 0, t<\lambda} \frac{\left(\frac{\|x+th^M\|_{l_p}^{p-2}}{\|x+th\|_{l_p}^{p-2}} - 1\right)x}{t}$$
$$= \frac{\|x^M\|_{l_p}^{p-2}}{\|x\|_{l_p}^{p-2}} h^M + \lim_{t\downarrow 0, t<\lambda} \frac{1}{\|x+th\|_{l_p}^{p-2}} \left(\frac{\|x+th^M\|_{l_p}^{p-2} - \|x+th\|_{l_p}^{p-2}}{t}\right)x$$
$$= \frac{\|x^M\|_{l_p}^{p-2}}{\|x\|_{l_p}^{p-2}} h^M + \frac{x}{\|x\|_{l_p}^{p-2}} \lim_{t\downarrow 0, t<\lambda} \frac{\|x+th^M\|_{l_p} - \|x+th\|_{l_p}}{t} \sum_{k=0}^{p-3} \|x+th^M\|_{l_p}^{p-3-k} \|x+th\|_{l_p}^{k}$$
$$= \frac{\|x^M\|_{l_p}^{p-2}}{\|x\|_{l_p}^{p-2}} h^M + \frac{x}{\|x\|_{l_p}^{p-2}} \sum_{k=0}^{p-3} \|x\|_{l_p}^{p-3-k} \|x\|_{l_p}^{k} \lim_{t\downarrow 0, t<\lambda} \frac{\|x+th^M\|_{l_p} - \|x+th\|_{l_p}}{t}$$
$$= \frac{\|x^M\|_{l_p}^{p-2}}{\|x\|_{l_p}^{p-2}} h^M + \frac{x}{\|x\|_{l_p}^{p-2}} (p-2)\|x\|_{l_p}^{p-3} \lim_{t\downarrow 0, t<\lambda} \frac{\|x+th^M\|_{l_p} - \|x+th\|_{l_p}}{t}$$
$$= \frac{\|x^M\|_{l_p}^{p-2}}{\|x\|_{l_p}^{p-2}} h^M + \frac{(p-2)x}{\|x\|_{l_p}} \lim_{t\downarrow 0, t<\lambda} \frac{\|x+th^M\|_{l_p} - \|x\|_{l_p}}{t}$$

$$= \frac{\|x^M\|_{l_p}^{p-2}}{\|x\|_{l_p}^{p-2}} h^M + \frac{(p-2)x}{\|x\|_{l_p}} \lim_{t\downarrow 0, t<\lambda} \frac{\|x+th^M\|_{l_p} - \|x\|_{l_p} - \left(\|x+th\|_{l_p} - \|x\|_{l_p}\right)}{t}$$
$$= \frac{\|x^M\|_{l_p}^{p-2}}{\|x\|_{l_p}^{p-2}} h^M + \frac{(p-2)}{\|x\|_{l_p}} \left(\Psi(x, h^M) - \Psi(x, h)\right)x.$$

Proof of part (i) in (b) for $h \in l_p^M$ with $h \neq \theta$. Let $x \in l_p^M \backslash B_M(r)$. Since $\|x\|_{l_p} > r$, there is $\delta > 0$, such that $x + th \in l_p^M \backslash B_M(r)$, for all $t \in (0, \delta)$. Noticing $x \neq \theta$, by part (b) in Theorem 6.2 and the definition of $\Psi$, we have

$$\Pi'_{B_M(r)}(x)(h)$$

$$
\begin{aligned}
&= \lim_{t\downarrow 0}\frac{\Pi_{B_M(r)}(x+th)-\Pi_{B_M(r)}(x)}{t}\\
&= \lim_{t\downarrow 0, t<\delta}\frac{\frac{r}{\|x+th\|_{l_p}}(x+th)-\frac{r}{\|x\|_{l_p}}x}{t}\\
&= \lim_{t\downarrow 0, t<\delta}\frac{\frac{rth}{\|x+th\|_{l_p}}+\left(\frac{r}{\|x+th\|_{l_p}}-\frac{r}{\|x\|_{l_p}}\right)x}{t}\\
&= \lim_{t\downarrow 0, t<\delta}\frac{rh}{\|x+th\|_{l_p}}+\lim_{t\downarrow 0, t<\delta}\frac{\left(\frac{r}{\|x+th\|_{l_p}}-\frac{r}{\|x\|_{l_p}}\right)x}{t}\\
&= \frac{r}{\|x\|_{l_p}}h - \lim_{t\downarrow 0, t<\delta}\frac{\frac{r}{\|x+th\|_{l_p}\|x\|_{l_p}}\left(\|x+th\|_{l_p}-\|x\|_{l_p}\right)x}{t}\\
&= \frac{r}{\|x\|_{l_p}}h - \frac{r}{\|x\|_{l_p}^2}\Psi(x,h)x.
\end{aligned}
$$

Proof of (b) for $h \in l_p\backslash l_p^M$. For any $x \in l_p^M\backslash B_M(r)$, it satisfies

$$
\frac{r}{\|x^M\|_{l_p}} = \frac{r}{\|x\|_{l_p}} < 1 = \frac{\|x\|_{l_p}^{p-2}}{\|x\|_{l_p}^{p-2}} = \frac{\|x^M\|_{l_p}^{p-2}}{\|x\|_{l_p}^{p-2}}. \tag{6.11}
$$

For any $h \in l_p\backslash l_p^M$(it is clear that $h \neq \theta$), by $r < \|x^M\|_{l_p} = \|x\|_{l_p}$ and by the continuity of the norm $\|\cdot\|_{l_p}$ and (6.11), there is $\lambda > 0$ such that $x + th \in l_p\backslash(C_M(r) \cup l_p^M)$ with

$$
r < \|(x+th)^M\|_{l_p} < \|x+th\|_{l_p},
$$

and

$$
\frac{\|(x+th)^M\|_{l_p}^{p-2}}{\|x+th\|_{l_p}^{p-2}} > \frac{r}{\|(x+th)^M\|_{l_p}}, \text{ for all } t \in (0, \lambda).
$$

Noticing $x \neq \theta$, by the first case of part (d) in Theorem 6.2 and (6.3), we have

$$
\begin{aligned}
&\Pi'_{B_M(r)}(x)(h)\\
&= \lim_{t\downarrow 0}\frac{\Pi_{B_M(r)}(x+th)-\Pi_{B_M(r)}(x)}{t}\\
&= \lim_{t\downarrow 0, t<\lambda}\frac{\frac{r}{\|(x+th)^M\|_{l_p}}(x+th)^M-\frac{r}{\|x\|_{l_p}}x}{t}\\
&= \lim_{t\downarrow 0, t<\lambda}\frac{\frac{rth^M}{\|(x+th)^M\|_{l_p}}+\left(\frac{r}{\|(x+th)^M\|_{l_p}}-\frac{r}{\|x\|_{l_p}}\right)x}{t}\\
&= \lim_{t\downarrow 0, t<\lambda}\frac{\frac{rth^M}{\|(x+th)^M\|_{l_p}}}{t}+\lim_{t\downarrow 0, t<\lambda}\frac{\left(\frac{r}{\|(x+th)^M\|_{l_p}}-\frac{r}{\|x\|_{l_p}}\right)x}{t}\\
&= \frac{r}{\|x\|_{l_p}}h^M - \frac{r}{\|x\|_{l_p}^2}\lim_{t\downarrow 0, t<\lambda}\frac{\|x+th^M\|_{l_p}-\|x\|_{l_p}}{t}x\\
&= \frac{r}{\|x\|_{l_p}}h^M - \frac{r}{\|x\|_{l_p}^2}\Psi(x,h^M)x.
\end{aligned}
$$

To complete the proof of part (b), for any $x \in l_p^M \backslash B_M(r)$, by Lemma 2.2, $\Psi(x, x) = \|x\|_{l_p}$. Then, from the fact that $\Psi(x, x) = \|x\|_{l_p}$ and taking $h = x$ in part (b) gets (ii) in (b).

Proof of (c). Let $x \in C_M^o(r) \backslash l_p^M$ satisfying $\|x^M\|_{l_p} < r$ and $\|x^M\|_{l_p} < \|x\|_{l_p}$. For $h \in l_p$ with $h \neq \theta$, there is $\delta > 0$, such that $x + th \in C_M^o(r) \backslash l_p^M$, for all $t \in (0, \delta)$. By part (c) in Theorem 6.2 and by (6.3), we have

$$
\begin{aligned}
&\Pi'_{B_M(r)}(x)(h) \\
&= \lim_{t \downarrow 0} \frac{\Pi_{B_M(r)}(x+th) - \Pi_{B_M(r)}(x)}{t} \\
&= \lim_{t \downarrow 0, t<\delta} \frac{(x+th)^M - x^M}{t} \\
&= \lim_{t \downarrow 0, t<\delta} \frac{th^M}{t} \\
&= h^M.
\end{aligned}
$$

Proof of (d). Let $x \in l_p \backslash (C_M(r) \cup l_p^M)$ with $\|x\|_{l_p} > \|x^M\|_{l_p} > r$. Proof of Part (i) in (d). Suppose $x$ satisfies

$$
\frac{\|x^M\|_{l_p}^{p-2}}{\|x\|_{l_p}^{p-2}} > \frac{r}{\|x^M\|_{l_p}}.
$$

In this case, for $h \in l_p$ with $h \neq \theta$, by the continuity of the norm $\|\cdot\|_{l_p}$, there is $\delta > 0$, such that

$$
\|x + th\|_{l_p} > \|(x + th)^M\|_{l_p} > r, \text{ for all } t \in (0, \delta).
$$

That is,

$$
x + th \in l_p \backslash (C_M(r) \cup l_p^M), \text{ for all } t \in (0, \delta),
$$

and

$$
\frac{\|(x+th)^M\|_{l_p}^{p-2}}{\|x+th\|_{l_p}^{p-2}} > \frac{r}{\|(x+th)^M\|_{l_p}}, \text{ for all } t \in (0, \delta).
$$

By the first case in part (d) in Theorem 6.2 and by (6.3), we have

$$
\begin{aligned}
&\Pi'_{B_M(r)}(x)(h) \\
&= \lim_{t \downarrow 0} \frac{\Pi_{B_M(r)}(x+th) - \Pi_{B_M(r)}(x)}{t} \\
&= \lim_{t \downarrow 0, t<\delta} \frac{\frac{r}{\|(x+th)^M\|_{l_p}}(x+th)^M - \frac{r}{\|x^M\|_{l_p}} x^M}{t} \\
&= \lim_{t \downarrow 0, t<\delta} \frac{\frac{rth^M}{\|(x+th)^M\|_{l_p}} + \left(\frac{r}{\|(x+th)^M\|_{l_p}} - \frac{r}{\|x^M\|_{l_p}}\right) x^M}{t} \\
&= \lim_{t \downarrow 0, t<\delta} \frac{rh^M}{\|(x+th)^M\|_{l_p}} + \lim_{t \downarrow 0, t<\delta} \frac{\left(\frac{r}{\|(x+th)^M\|_{l_p}} - \frac{r}{\|x^M\|_{l_p}}\right) x^M}{t} \\
&= \frac{r}{\|x^M\|_{l_p}} h^M - \frac{r}{\|x^M\|_{l_p}^2} \lim_{t \downarrow 0, t<\delta} \frac{\|x^M + th^M\|_{l_p} - \|x^M\|_{l_p}}{t} x^M \\
&= \frac{r}{\|x^M\|_{l_p}} h^M - \frac{r}{\|x^M\|_{l_p}^2} \Psi(x^M, h^M) x^M.
\end{aligned}
$$

Proof of part (ii) in (d). For $x \in l_p \backslash (C_M(r) \cup l_p^M)$ with $\|x\|_{l_p} > \|x^M\|_{l_p} > r$, suppose $x$ satisfies

$$\frac{\|x^M\|_{l_p}^{p-2}}{\|x\|_{l_p}^{p-2}} < \frac{r}{\|x^M\|_{l_p}}.$$

In this case, for $h \in l_p$ with $h \neq \theta$, by the continuity of the norm $\|\cdot\|_{l_p}$, there is $\lambda > 0$, such that

$$x + th \in l_p \backslash (C_M(r) \cup l_p^M), \text{ for all } t \in (0, \lambda),$$

and

$$\frac{\|(x+th)^M\|_{l_p}^{p-2}}{\|x+th\|_{l_p}^{p-2}} < \frac{r}{\|(x+th)^M\|_{l_p}}, \text{ for all } t \in (0, \lambda).$$

By the second part in (d) in Theorem 6.2 and by (6.3), for $p$ to be an integer greater than 1, we have

$$\begin{aligned}
&\Pi'_{B_M(r)}(x)(h) \\
&= \lim_{t\downarrow 0} \frac{\Pi_{B_M(r)}(x+th) - \Pi_{B_M(r)}(x)}{t} \\
&= \lim_{t\downarrow 0, t<\lambda} \frac{\frac{\|(x+th)^M\|_{l_p}^{p-2}}{\|x+th\|_{l_p}^{p-2}}(x+th)^M - \frac{\|x^M\|_{l_p}^{p-2}}{\|x\|_{l_p}^{p-2}} x^M}{t} \\
&= \lim_{t\downarrow 0, t<\lambda} \frac{\frac{\|(x+th)^M\|_{l_p}^{p-2}}{\|x+th\|_{l_p}^{p-2}} th^M + \left(\frac{\|(x+th)^M\|_{l_p}^{p-2}}{\|x+th\|_{l_p}^{p-2}} - \frac{\|x^M\|_{l_p}^{p-2}}{\|x\|_{l_p}^{p-2}}\right) x^M}{t} \\
&= \frac{\|x^M\|_{l_p}^{p-2}}{\|x\|_{l_p}^{p-2}} h^M + \lim_{t\downarrow 0, t<\lambda} \frac{\left(\frac{\|(x+th)^M\|_{l_p}^{p-2}}{\|x+th\|_{l_p}^{p-2}} - \frac{\|x^M\|_{l_p}^{p-2}}{\|x\|_{l_p}^{p-2}}\right) x^M}{t} \\
&= \frac{\|x^M\|_{l_p}^{p-2}}{\|x\|_{l_p}^{p-2}} h^M + \frac{x^M}{\|x\|_{l_p}^{2(p-2)}} \lim_{t\downarrow 0, t<\lambda} \frac{\|x^M + th^M\|_{l_p}^{p-2} \|x\|_{l_p}^{p-2} - \|x+th\|_{l_p}^{p-2} \|x^M\|_{l_p}^{p-2}}{t} \\
&= \frac{\|x^M\|_{l_p}^{p-2}}{\|x\|_{l_p}^{p-2}} h^M + \frac{x^M}{\|x\|_{l_p}^{2(p-2)}} \lim_{t\downarrow 0, t<\lambda} \frac{\|x\|_{l_p}^{p-2}\left(\|x^M + th^M\|_{l_p}^{p-2} - \|x^M\|_{l_p}^{p-2}\right) - \|x^M\|_{l_p}^{p-2}\left(\|x+th\|_{l_p}^{p-2} - \|x\|_{l_p}^{p-2}\right)}{t} \\
&= \frac{\|x^M\|_{l_p}^{p-2}}{\|x\|_{l_p}^{p-2}} h^M + \frac{x^M}{\|x\|_{l_p}^{2(p-2)}} \left(\|x\|_{l_p}^{p-2} \lim_{t\downarrow 0, t<\lambda} \frac{\|x^M + th^M\|_{l_p}^{p-2} - \|x^M\|_{l_p}^{p-2}}{t} - \|x^M\|_{l_p}^{p-2} \lim_{t\downarrow 0, t<\lambda} \frac{\|x+th\|_{l_p}^{p-2} - \|x\|_{l_p}^{p-2}}{t}\right) \\
&= \frac{\|x^M\|_{l_p}^{p-2}}{\|x\|_{l_p}^{p-2}} h^M + \frac{x^M}{\|x\|_{l_p}^{2(p-2)}} \left(\|x\|_{l_p}^{p-2} \lim_{t\downarrow 0, t<\lambda} \frac{\|x^M + th^M\|_{l_p} - \|x^M\|_{l_p}}{t} \sum_{k=0}^{p-3} \|x^M + th^M\|_{l_p}^{p-3-k} \|x^M\|_{l_p}^{k}\right. \\
&\qquad \left. - \|x^M\|_{l_p}^{p-2} \lim_{t\downarrow 0, t<\lambda} \frac{\|x+th\|_{l_p} - \|x\|_{l_p}}{t} \sum_{k=0}^{p-3} \|x + th\|_{l_p}^{p-3-k} \|x\|_{l_p}^{k}\right) \\
&= \frac{\|x^M\|_{l_p}^{p-2}}{\|x\|_{l_p}^{p-2}} h^M + \frac{1}{\|x\|_{l_p}^{2(p-2)}} \left(\|x\|_{l_p}^{p-2}(p-2)\|x^M\|_{l_p}^{p-3}\Psi(x^M, h^M) - \|x^M\|_{l_p}^{p-2}(p-2)\|x\|_{l_p}^{p-3}\Psi(x, h)\right) x^M \\
&= \frac{\|x^M\|_{l_p}^{p-2}}{\|x\|_{l_p}^{p-2}} h^M + \frac{p-2}{\|x\|_{l_p}^{2(p-2)}} \left(\|x\|_{l_p}^{p-2}\|x^M\|_{l_p}^{p-3}\Psi(x^M, h^M) - \|x^M\|_{l_p}^{p-2}\|x\|_{l_p}^{p-3}\Psi(x, h)\right) x^M
\end{aligned}$$

$$= \frac{\|x^M\|_{l_p}^{p-2}}{\|x\|_{l_p}^{p-2}} h^M + \frac{(p-2)\|x^M\|_{l_p}^{p-2}}{\|x\|_{l_p}^{p-2}} \left( \frac{\Psi(x^M, h^M)}{\|x^M\|_{l_p}} - \frac{\Psi(x,h)}{\|x\|_{l_p}} \right) x^M. \qquad \square$$

For a given nonempty subset $M$ of $\mathbb{N}$, if $\mathbb{N} \backslash M = \emptyset$, then $l_p^M$ coincides with $l_p$. In this case, we have

$$B_M(r) = C_M(r) = B(r), \text{ and } h^M = h, \text{ for any } h \in l_p.$$

Then, by Theorem 6.5, obtain

**Corollary 6.6**. *For any* $r > 0$ *and for any* $x, h \in l_p$ *with* $h \neq \theta$, *we have*

(a) *For* $x \in B^o(r)$,

$$\Pi'_{B(r)}(x)(h) = h;$$

(b) *For* $x \in l_p \backslash B(r)$,

$$\Pi'_{B(r)}(x)(h) = \frac{r}{\|x\|_{l_p}} h - \frac{r}{\|x\|_{l_p}^2} \Psi(x, h) x.$$

*In particular,*

$$\Pi'_{B(r)}(x)(x) = \theta, \text{ for any } x \in l_p \backslash B(r).$$

### 6.3. The standard metric projection *P* onto closed balls in $l_p$

In the previous two subsections, we calculate the solutions of the generalized metric projection $\Pi_{B_M(r)}$ in $\boldsymbol{l_p}$. To find the connections of the directional differentiability between the generalized metric projection $\Pi$ and the (standard) metric projection $P$ from $\boldsymbol{l_p}$ onto a closed ball $B_M(r)$ in $l_p^M$, in this subsection, we will calculate the solutions of $P_{B_M(r)}$ in $\boldsymbol{l_p}$. Recall that the basic variational principle of the (standard) metric projection $P_{B_M(r)}$ is: for any $x \in \boldsymbol{l_p}$ and $u \in B_M(r)$,

$$u = P_{B_M(r)}(x) \quad \Leftrightarrow \quad \langle J(x-u), u-y \rangle \geq 0, \text{ for all } y \in B_M(r). \tag{6.12}$$

**Theorem 6.7**. *For any* $r > 0$ *and* $x \in \boldsymbol{l_p}$, *we have*

(a) $P_{B_M(r)}(x) = x$, *for* $x \in B_M(r)$;
(b) $P_{B_M(r)}(x) = \frac{r}{\|x\|_{l_p}} x$, *for* $x \in l_p^M \backslash B_M(r)$;
(c) $P_{B_M(r)}(x) = x^M$, *for* $x \in C_M(r) \backslash l_p^M$;
(d) $P_{B_M(r)}(x) = \frac{r}{\|x^M\|_{l_p}} x^M$, *for* $x \in l_p \backslash (C_M(r) \cup l_p^M)$.

*Proof*. Part (a) is clear. We prove part (b). For any given $x \in l_p^M \backslash B_M(r)$, $x$ must satisfies $\|x^M\|_{l_p} = \|x\|_{l_p} > r$. Then, for any $y \in B_M(r)$, we have

$$\begin{aligned}
&\left\langle J\left(x - \frac{r}{\|x\|_{l_p}} x\right), \frac{r}{\|x\|_{l_p}} x - y \right\rangle \\
&= \left(1 - \frac{r}{\|x\|_{l_p}}\right) \left\langle J(x), \frac{r}{\|x\|_{l_p}} x - y \right\rangle \\
&= \left(1 - \frac{r}{\|x\|_{l_p}}\right) \left(r\|x\|_{l_p} - \langle J(x), y \rangle\right) \\
&\geq \left(\|x\|_{l_p} - r\right) \left(r - \|y\|_{l_p}\right) \\
&\geq 0, \text{ for all } y \in B_M(r).
\end{aligned}$$

By $x \in l_p^M$ and $\|x\|_{l_p} > r$, we have $\frac{r}{\|x\|_{l_p}} x \in S_M(r)$. By the basic variational principle of $P_{B_M(r)}$ in (6.12), the above inequality implies

$$P_{B_M(r)}(x) = \frac{r}{\|x\|_{l_p}} x, \quad \text{for any } x \in l_p^M \backslash B_M(r).$$

Proof of (c). For $x \in C_M(r) \backslash l_p^M$, we have

$$x^M \in B_M(r) \text{ and } x \notin l_p^M. \tag{6.13}$$

It follows that

$$0 < \|x^M\|_{l_p} \le r \quad \text{and} \quad \|x^M\|_{l_p} < \|x\|_{l_p}.$$

Then, for any $y \in B_M(r) \subseteq l_p^M$, we have $x^M - y \in l_p^M$. By (6.2) and Lemma 6.1, we obtain

$$\begin{aligned} &\langle J(x - x^M),\ x^M - y\rangle \\ &= \langle J(x^{\mathbb{N}\backslash M}), x^M - y\rangle \\ &= 0, \text{ for all } y \in B_M(r). \end{aligned}$$

By (6.13) and the basic variational principle of $P_{B_M(r)}$, this implies $P_{B_M(r)}(x) = x_M$. Then, we prove part (d). For $x \in l_p \backslash (C_M(r) \cup l_p^M)$, we have $x^M \notin B_M(r)$ and $x \notin l_p^M$. It follows that

$$r < \|x^M\|_{l_p} < \|x\|_{l_p}. \tag{6.14}$$

Then, for any $y \in B_M(r) \subseteq l_p^M$, by (6.1) and (6.14), we have

$$\begin{aligned} &\langle J\left(x - \frac{r}{\|x^M\|_{l_p}} x^M\right),\ \frac{r}{\|x^M\|_{l_p}} x^M - y\rangle \\ &= \langle J\left(\left(1 - \frac{r}{\|x^M\|_{l_p}}\right) x^M + x^{\mathbb{N}\backslash M}\right),\ \frac{r}{\|x^M\|_{l_p}} x^M - y\rangle \end{aligned}$$

Let $a(x) = \left(1 - \frac{r}{\|x^M\|_{l_p}}\right) x^M + x^{\mathbb{N}\backslash M}$. It is clear that $a(x) \in l_p \backslash \{\theta\}$. By $x^M - y \in l_p^M$, by Lemma 6.1, and $1 - \frac{r}{\|x^M\|_{l_p}} > 0$, we have

$$\begin{aligned} &\langle J\left(x - \frac{r}{\|x^M\|_{l_p}} x^M\right),\ \frac{r}{\|x^M\|_{l_p}} x^M - y\rangle \\ &= \langle \frac{\left\|\left(1 - \frac{r}{\|x^M\|_{l_p}}\right) x^M\right\|_{l_p}^{p-2}}{\|a(x)\|_{l_p}^{p-2}} J\left(\left(1 - \frac{r}{\|x^M\|_{l_p}}\right) x^M\right) + \frac{\|x^{\mathbb{N}\backslash M}\|_{l_p}^{p-2}}{\|a(x)\|_{l_p}^{p-2}} J\left(x^{\mathbb{N}\backslash M}\right),\ \frac{r}{\|x^M\|_{l_p}} x^M - y\rangle \\ &= \langle \frac{\left(1 - \frac{r}{\|x^M\|_{l_p}}\right)\left\|\left(1 - \frac{r}{\|x^M\|_{l_p}}\right) x^M\right\|_{l_p}^{p-2}}{\|a(x)\|_{l_p}^{p-2}} J(x^M) + \frac{\|x^{\mathbb{N}\backslash M}\|_{l_p}^{p-2}}{\|a(x)\|_{l_p}^{p-2}} J\left(x^{\mathbb{N}\backslash M}\right),\ \frac{r}{\|x^M\|_{l_p}} x^M - y\rangle \end{aligned}$$

$$= \left\langle \frac{\left(1-\frac{r}{\|x^M\|_{l_p}}\right)\left\|\left(1-\frac{r}{\|x^M\|_{l_p}}\right)x^M\right\|_{l_p}^{p-2}}{\|a(x)\|_{l_p}^{p-2}} J(x^M),\ \frac{r}{\|x^M\|_{l_p}} x^M - y \right\rangle$$

$$= \frac{\left(1-\frac{r}{\|x^M\|_{l_p}}\right)^{p-1} \|x^M\|_{l_p}^{p-2}}{\|a(x)\|_{l_p}^{p-2}} \left\langle J(x^M),\ \frac{r}{\|x^M\|_{l_p}} x^M - y \right\rangle$$

$$= \frac{\left(1-\frac{r}{\|x^M\|_{l_p}}\right)^{p-1} \|x^M\|_{l_p}^{p-2}}{\|a(x)\|_{l_p}^{p-2}} \left(r\|x^M\|_{l_p} - \langle J(x^M), y\rangle\right)$$

$$\geq \frac{\left(1-\frac{r}{\|x^M\|_{l_p}}\right)^{p-1} \|x^M\|_{l_p}^{p-2}}{\|a(x)\|_{l_p}^{p-2}} \|x^M\|_{l_p} \left(r - \|y\|_{l_p}\right)$$

$\geq 0$, for all $y \in B_M(r)$.

By $\|x^M\|_{l_p} > r$, we have $\frac{r}{\|x^M\|_{l_p}} x^M \in S_M(r)$. By the basic variational principle of $P_{B_M(r)}$, this implies

$$P_{B_M(r)}(x) = \frac{r}{\|x^M\|_{l_p}} x^M, \quad \text{for } x \in l_p \backslash (C_M(r) \cup l_p^M).$$ □

Notice that if $\mathbb{N}\backslash M = \emptyset$, then $l_p^M$ coincides with $l_p$. In this case, we have $B_M(r) = C_M(r) = B(r)$. Then, as a consequence of Theorem 6.7 and Corollary 6.3, we have

**Corollary 6.8**. *For any $r > 0$ and $x \in l_p$, we have*

(a) $P_{B(r)}(x) = \Pi_{B(r)}(x) = x$, *for* $x \in B(r)$;
(b) $P_{B(r)}(x) = \Pi_{B(r)}(x) = \frac{r}{\|x\|_{l_p}} x$, *for* $x \in l_p \backslash B(r)$.

**Remarks 6.9**. By respectively comparing parts (c) and (d) in Theorem 6.2 and parts (c) and (d) in Theorem 6.7, we see the big differences between $\Pi_{B_M(r)}$ and $P_{B_M(r)}$.

(i) If $\mathbb{N}\backslash M \neq \emptyset$, then, in general,

$$P_{B_M(r)}(x) \neq \Pi_{B_M(r)}(x), \quad \text{for } x \in C_M(r) \backslash l_p^M \text{ or } x \in l_p \backslash (C_M(r) \cup l_p^M).$$

(ii) If $\mathbb{N}\backslash M = \emptyset$, by comparing parts (a) and (b) in Corollary 6.8, and parts (a) and (b) in Corollary 6.4. respectively, we have

$$P_{B(r)}(x) = \Pi_{B(r)}(x), \quad \text{for any } x \in l_p.$$

Now, by Theorem 6.7, we study the directional differentiability of $P_{B_M(r)}$.

**Theorem 6.10**. *For any $r > 0$ and $x, h \in l_p$ with $h \neq \theta$, we have,*

(a) *If $x \in B_M^o(r)$, then*

$$P'_{B_M(r)}(x)(h) = \begin{cases} h, & \text{for } h \in l_p^M \\ h^M, & \text{for } h \notin l_p^M \end{cases};$$

(b) *If* $x \in l_p^M \backslash B_M(r)$, *then*

$$P'_{B_M(r)}(x)(h) = \begin{cases} \frac{r}{\|x\|_{l_p}}\left(h - \frac{\Psi(x,h)}{\|x\|_{l_p}}x\right), & \text{for } h \in l_p^M \\ \frac{r}{\|x\|_{l_p}}\left(h^M - \frac{\Psi(x,h^M)}{\|x\|_{l_p}}x\right), & \text{for } h \notin l_p^M \end{cases}.$$

*In particular,*

$$P'_{B_M(r)}(x)(x) = \theta, \textit{ for any } x \in l_p^M \backslash B_M(r);$$

(c) *If* $x \in C_M^o(r) \backslash l_p^M$, *then*

$$P'_{B_M(r)}(x)(h) = h^M;$$

(d) *If* $x \in l_p \backslash (C_M(r) \cup l_p^M)$, *then*

$$P'_{B_M(r)}(x)(h) = \frac{r}{\|x^M\|_{l_p}}\left(h^M - \frac{\Psi(x^M, h^M)}{\|x^M\|_{l_p}}x^M\right).$$

*In particular,*

$$P'_{B_M(r)}(x)(x) = \theta, \textit{ for any } x \in l_p \backslash (C_M(r) \cup l_p^M).$$

*Proof*. Proof of (a). Suppose $x \in B_M^o(r)$. If $h \in l_p^M$, then, there is $\delta > 0$, such that $x + th \in B_M(r)$, for all $t \in (0, \delta)$. By part (a) in Theorem 6.7, we have

$$\begin{aligned} & P'_{B_M(r)}(x)(h) \\ & = \lim_{t \downarrow 0} \frac{P_{B_M(r)}(x+th) - P_{B_M(r)}(x)}{t} \\ & = \lim_{t \downarrow 0, t<\delta} \frac{x+th-x}{t} \\ & = \lim_{t \downarrow 0, t<\delta} \frac{th}{t} \\ & = h. \end{aligned}$$

If $h \notin l_p^M$, that is, $h \in l_p \backslash l_p^M$, then, there is $\lambda > 0$, such that $x + th \in C_M(r) \backslash l_p^M$, for all $t \in (0, \lambda)$. By parts (a, c) in Theorem 6.7 and by Lamma 6.1, we have

$$\begin{aligned} & P'_{B_M(r)}(x)(h) \\ & = \lim_{t \downarrow 0} \frac{P_{B_M(r)}(x+th) - P_{B_M(r)}(x)}{t} \\ & = \lim_{t \downarrow 0, t<\lambda} \frac{(x+th)^M - x}{t} \\ & = \lim_{t \downarrow 0, t<\lambda} \frac{x+th^M - x}{t} \\ & = \lim_{t \downarrow 0, t<\lambda} \frac{th^M}{t} \\ & = h^M. \end{aligned}$$

Proof of (b). Let $x \in l_p^M \backslash B_M(r)$. It implies $\|x^M\|_{l_p} = \|x\|_{l_p} > r$. If $h \in l_p^M$ with $h \neq \theta$, then, there is $\delta > 0$, such that $x + th \in l_p^M \backslash B_M(r)$, for all $t \in (0, \delta)$. Noticing $x \neq \theta$, by part (b) in Theorem 6.7 and the definition of $\Psi$, we have

$$\begin{aligned} & P'_{B_M(r)}(x)(h) \\ & = \lim_{t \downarrow 0} \frac{P_{B_M(r)}(x+th) - P_{B_M(r)}(x)}{t} \end{aligned}$$

$$= \lim_{t\downarrow 0, t<\delta} \frac{\frac{r}{\|x+th\|_{l_p}}(x+th) - \frac{r}{\|x\|_{l_p}}x}{t}$$
$$= \lim_{t\downarrow 0, t<\delta} \frac{\frac{rth}{\|x+th\|_{l_p}} + \left(\frac{r}{\|x+th\|_{l_p}} - \frac{r}{\|x\|_{l_p}}\right)x}{t}$$
$$= \lim_{t\downarrow 0, t<\delta} \frac{rh}{\|x+th\|_{l_p}} + \lim_{t\downarrow 0, t<\delta} \frac{\left(\frac{r}{\|x+th\|_{l_p}} - \frac{r}{\|x\|_{l_p}}\right)x}{t}$$
$$= \frac{r}{\|x\|_{l_p}}h - \lim_{t\downarrow 0, t<\delta} \frac{\frac{r}{\|x+th\|_{l_p}\|x\|_{l_p}}\left(\|x+th\|_{l_p} - \|x+th\|_{l_p}\right)x}{t}$$
$$= \frac{r}{\|x\|_{l_p}}h - \frac{r}{\|x\|_{l_p}^2}\Psi(x,h)x$$
$$= \frac{r}{\|x\|_{l_p}}\left(h - \frac{\Psi(x,h)}{\|x\|_{l_p}}x\right).$$

For any given $x \in l_p^M \backslash B_M(r)$ with $\|x^M\|_{l_p} = \|x\|_{l_p} > r$, if $h \in l_p \backslash l_p^M$, then, there is $\lambda > 0$, such that $x + th \in l_p \backslash (C_M(r) \cup l_p^M)$, for all $t \in (0, \lambda)$. Then, by parts (b, d) in Theorem 6.7 and Lamma 6.1, we have

$$P'_{B_M(r)}(x)(h)$$
$$= \lim_{t\downarrow 0} \frac{P_{B_M(r)}(x+th) - P_{B_M(r)}(x)}{t}$$
$$= \lim_{t\downarrow 0, t<\lambda} \frac{\frac{r}{\|(x+th)^M\|_{l_p}}(x+th)^M - \frac{r}{\|x\|_{l_p}}x}{t}$$
$$= \lim_{t\downarrow 0, t<\lambda} \frac{\frac{rth^M}{\|(x+th)^M\|_{l_p}} + \left(\frac{r}{\|(x+th)^M\|_{l_p}} - \frac{r}{\|x\|_{l_p}}\right)x}{t}$$
$$= \lim_{t\downarrow 0, t<\lambda} \frac{rh_M}{\|(x+th)^M\|_{l_p}} + \lim_{t\downarrow 0, t<\lambda} \frac{\left(\frac{r}{\|(x+th)^M\|_{l_p}} - \frac{r}{\|x\|_{l_p}}\right)x}{t}$$
$$= \frac{r}{\|x\|_{l_p}}h^M - \frac{r}{\|x\|_{l_p}^2}\Psi(x,h^M)x$$
$$= \frac{r}{\|x\|_{l_p}}\left(h^M - \frac{\Psi(x,h^M)}{\|x\|_{l_p}}x\right).$$

Proof of (c). Let $x \in C_M^o(r) \backslash l_p^M$ with $\|x^M\|_{l_p} < r$ and $\|x^M\|_{l_p} < \|x\|_{l_p}$. For $h \in l_p$ with $h \neq \theta$, there is $\delta > 0$, such that $x + th \in C_M^o(r) \backslash l_p^M$, for all $t \in (0, \delta)$. By part (c) in Theorem 6.7 and by Lamma 6.1, we have

$$P'_{B_M(r)}(x)(h)$$
$$= \lim_{t\downarrow 0} \frac{P_{B_M(r)}(x+th) - P_{B_M(r)}(x)}{t}$$
$$= \lim_{t\downarrow 0, t<\delta} \frac{(x+th)^M - x^M}{t}$$
$$= \lim_{t\downarrow 0, t<\delta} \frac{th^M}{t}$$
$$= h^M.$$

Proof of (d). Let $x \in l_p \backslash (C_M(r) \cup l_p^M)$ with $\|x\|_{l_p} > \|x^M\|_{l_p} > r$. For $h \in l_p$ with $h \neq \theta$, there is $\delta > 0$, such that $\|(x + th)\|_{l_p} > \|(x + th)^M\|_{l_p} > r$, that is, $x + th \in l_p \backslash (C_M(r) \cup l_p^M)$, for all $t \in (0, \delta)$. By part (d) in Theorem 6.7 and by Lamma 6.1, we have

$$
\begin{aligned}
&P'_{B_M(r)}(x)(h)\\
&= \lim_{t\downarrow 0} \frac{P_{B_M(r)}(x+th)-P_{B_M(r)}(x)}{t}\\
&= \lim_{t\downarrow 0, t<\lambda} \frac{\frac{r}{\|(x+th)^M\|_{l_p}}(x+th)^M - \frac{r}{\|x^M\|_{l_p}}x^M}{t}\\
&= \lim_{t\downarrow 0, t<\lambda} \frac{\frac{rth^M}{\|(x+th)^M\|_{l_p}} + \left(\frac{r}{\|(x+th)^M\|_{l_p}} - \frac{r}{\|x^M\|_{l_p}}\right)x^M}{t}\\
&= \lim_{t\downarrow 0, t<\lambda} \frac{rh_M}{\|(x+th)^M\|_{l_p}} + \lim_{t\downarrow 0, t<\lambda} \frac{\left(\frac{r}{\|(x+th)^M\|_{l_p}} - \frac{r}{\|x^M\|_{l_p}}\right)x^M}{t}\\
&= \frac{r}{\|x^M\|_{l_p}} h^M - \frac{r}{\|x^M\|_{l_p}^2} \Psi(x^M, h^M) x^M\\
&= \frac{r}{\|x^M\|_{l_p}} \left(h^M - \frac{\Psi(x^M, h^M)}{\|x^M\|_{l_p}} x^M\right). \qquad \square
\end{aligned}
$$

In particular, if $\mathbb{N}\backslash M = \emptyset$, then, $B_M(r) = C_M(r) = B(r)$ and $x^M = x$, for any $x \in H$. As a corollary of Theorem 6.7 and Corollary 6.6, we have

**Corollary 6.11**. *Let* $r > 0$. *For* $x, h \in l_p$ *with* $h \neq \theta$, *we have,*

(a) *If* $x \in B^o(r)$, *then*

$$P'_{B(r)}(x)(h) = \Pi'_{B(r)}(x)(h) = h;$$

(b) *If* $x \in l_p\backslash B(r)$, *then*

$$P'_{B(r)}(x)(h) = \Pi'_{B(r)}(x)(h) = \frac{r}{\|x\|_{l_p}}\left(h - \frac{\Psi(x,h)}{\|x\|_{l_p}} x\right).$$

*In particular,*

$$P'_{B(r)}(x)(x) = \Pi'_{B(r)}(x)(h) = \theta, \textit{ for any } x \in l_p\backslash B(r).$$

**Remarks 6.12**. By correspondingly comparing parts (c) and (d) in Theorem 6.5 and Theorem 6.10, there are differences between $\Pi'_{B_M(r)}$ and $P'_{B_M(r)}$. That is,

(i) If $\mathbb{N}\backslash M \neq \emptyset$, then, in general,

$$P'_{B_M(r)} \neq \Pi'_{B_M(r)}, \text{ for } x \in l_p.$$

(ii) If $\mathbb{N}\backslash M = \emptyset$, by comparing parts (a) and (b) in Corollary 6.11, respectively, we have

$$P'_{B(r)} = \Pi'_{B(r)}, \text{ for any } x \in l_p.$$

### 6.4. The generalized metric projection Π onto closed and convex cylinders in $l_3$

In this subsection, we study the generalized metric projection Π onto closed and convex cylinders of $\boldsymbol{l_p}$ with respect to a special case for $p = 3$. Let $B_M(1)$ be the closed unit ball in $l_3^M$ and let $C_M(1)$ be the closed and convex cylinder in $l_3$ with base $B_M(1)$ in $l_3^M$. If $\mathbb{N}\backslash M \neq \emptyset$, then $C_M(1)$ is a closed, unbounded and convex subset in $l_3$. In this subsection, we calculate the values of the generalized metric projection $\Pi_{C_M(1)}$ from $l_3$ onto $C_M(1)$. In contrast with $\Pi_{B_M(1)}$ studied in subsection 6.1, we will see that, even in the special case of $p = 3$, the solutions of $\Pi_{C_M(1)}$ are more complicated than the solutions of $P_{B_M(1)}$. We need the following lemmas.

**Lemma 6.13**. *For any* $x \in l_3$ *with* $\|x^M\|_{l_3} > 1$ *and* $\left\|x^{\mathbb{N}\setminus M}\right\|_{l_3} > 0$. *Let* $b = b(x)$ *and* $a = a(x)$ *defined by*

$$b = \left(\frac{\|x^{\mathbb{N}\setminus M}\|_{l_3}^3 + \sqrt{\|x^{\mathbb{N}\setminus M}\|_{l_3}^6 + 4\|x\|_{l_3}^3}}{2\|x\|_{l_3}^3}\right)^{\frac{1}{3}}. \tag{6.15}$$

*and*

$$a = a(x) = \frac{1}{\|x^M\|_{l_3}} x^M + b x^{\mathbb{N}\setminus M}.$$

*Then,* $a(x)$ *satisfies*

(a) $a^M \in S_M(r)$;
(b) $(J(x))^{\mathbb{N}\setminus M} = (J(a))^{\mathbb{N}\setminus M}$.

Proof. We only prove part (b). By (6.3) with $p - 2 = 1$, we have

$$J(x) = \frac{\|x^M\|_{l_3}}{\|x\|_{l_3}} J(x^M) + \frac{\|x^{\mathbb{N}\setminus M}\|_{l_3}}{\|x\|_{l_3}} J(x^{\mathbb{N}\setminus M}), \tag{6.3}$$

and

$$\|x\|_{l_3}^3 = \|x^M\|_{l_3}^3 + \left\|x^{\mathbb{N}\setminus M}\right\|_{l_3}^3.$$

For a given $b > 0$, let

$$a = a(x) = \frac{1}{\|x^M\|_{l_3}} x^M + b x^{\mathbb{N}\setminus M}.$$

Then

$$\|a\|_{l_3} = \left(\left(\frac{1}{\|x^M\|_{l_3}}\right)^3 \|x^M\|_{l_3}^3 + b^3 \left\|x^{\mathbb{N}\setminus M}\right\|_{l_3}^3\right)^{\frac{1}{3}}$$

$$= \left(1 + b^3 \left\|x^{\mathbb{N}\setminus M}\right\|_{l_3}^3\right)^{\frac{1}{3}}. \tag{6.16}$$

By (6.3), we have

$$J(a) = \frac{\left\|\frac{1}{\|x^M\|_{l_3}} x^M\right\|_{l_3}^{3-2}}{\|a\|_{l_3}^{3-2}} J\left(\frac{1}{\|x^M\|_{l_3}} x^M\right) + \frac{\|b x^{\mathbb{N}\setminus M}\|_{l_3}^{3-2}}{\|a\|_{l_3}^{3-2}} J(b x^{\mathbb{N}\setminus M})$$

$$= \frac{1}{\|a\|_{l_3}} \frac{1}{\|x^M\|_{l_3}} J(x^M) + \frac{b\|x^{\mathbb{N}\setminus M}\|_{l_3}}{\|a\|_{l_3}} b J(x^{\mathbb{N}\setminus M})$$

$$= \frac{1}{\|x^M\|_{l_3}\|a\|_{l_3}} J(x^M) + \frac{b^2\|x^{\mathbb{N}\setminus M}\|_{l_3}}{\|a\|_{l_3}} J(x^{\mathbb{N}\setminus M}). \tag{6.17}$$

In order to have $(J(x))^{\mathbb{N}\setminus M} = (J(a))^{\mathbb{N}\setminus M}$, by (6.3) and (6.17), we need

$$\frac{\|x^{\mathbb{N}\setminus M}\|_{l_3}}{\|x\|_{l_3}} = \frac{b^2\|x^{\mathbb{N}\setminus M}\|_{l_3}}{\|a\|_{l_3}}.$$

That is,

$$\frac{1}{\|x\|_{l_3}} = \frac{b^2}{\left(r^3+b^3\|x^{\mathbb{N}\backslash M}\|_{l_3}^3\right)^{\frac{1}{3}}}.$$

It follows that

$$b^6\|x\|_{l_3}^3 - b^3\left\|x^{\mathbb{N}\backslash M}\right\|_{l_3}^3 - 1 = 0.$$

Since we need $b > 0$, we obtain

$$b^3 = \frac{\|x^{\mathbb{N}\backslash M}\|_{l_3}^3+\sqrt{\|x^{\mathbb{N}\backslash M}\|_{l_3}^6+4\|x\|_{l_3}^3}}{2\|x\|_{l_3}^3}. \qquad \square$$

**Lemma 6.14**. *For any* $x \in l_3$ *with* $\|x^M\|_{l_3} > 1$ *and* $\left\|x^{\mathbb{N}\backslash M}\right\|_{l_3} > 0$, *let* $b = b(x)$ *and* $a = a(x)$ *be defined as in Lemma* 6.13. *If x satisfies*

$$\|x\|_{l_3} \le \left(\frac{1+\sqrt{5}}{2}\right)^{\frac{1}{3}} \left\|x^{\mathbb{N}\backslash M}\right\|_{l_3}^2, \tag{6.18}$$

*then*,

$$\frac{\|x^M\|_{l_3}}{\|x\|_{l_3}} > \frac{1}{\|x^M\|_{l_3}\|a\|_{l_3}}.$$

*Proof*. By (6.15) and (6.16) in Lemma 6.13, we have

$$\begin{aligned}
&\|a\|_{l_3}^3 \\
&= 1 + b^3\left\|x^{\mathbb{N}\backslash M}\right\|_{l_3}^3 \\
&= 1 + \frac{\|x^{\mathbb{N}\backslash M}\|_{l_3}^3+\sqrt{\|x^{\mathbb{N}\backslash M}\|_{l_3}^6+4\|x\|_{l_3}^3}}{2\|x\|_{l_3}^3}\left\|x^{\mathbb{N}\backslash M}\right\|_{l_3}^3
\end{aligned}$$

Then, we calculate

$$\begin{aligned}
&\frac{\|x^M\|_{l_3}}{\|x\|_{l_3}} > \frac{1}{\|x^M\|_{l_3}\|a\|_{l_3}} \\
\Leftrightarrow\quad & \|x^M\|_{l_3}^6\,\|a\|_{l_3}^3 > \|x\|_{l_3}^3 \\
\Leftrightarrow\quad & \|x^M\|_{l_3}^6\left(1 + \frac{\|x^{\mathbb{N}\backslash M}\|_{l_3}^3+\sqrt{\|x^{\mathbb{N}\backslash M}\|_{l_3}^6+4\|x\|_{l_3}^3}}{2\|x\|_{l_3}^3}\left\|x^{\mathbb{N}\backslash M}\right\|_{l_3}^3\right) > \|x\|_{l_3}^3 \\
\Leftrightarrow\quad & \|x^M\|_{l_3}^6\left(2\|x\|_{l_3}^3 + \left(\left\|x^{\mathbb{N}\backslash M}\right\|_{l_3}^3 + \sqrt{\|x^{\mathbb{N}\backslash M}\|_{l_3}^6 + 4\|x\|_{l_3}^3}\right)\left\|x^{\mathbb{N}\backslash M}\right\|_{l_3}^3\right) > 2\|x\|_{l_3}^6 \\
\Leftrightarrow\quad & \left(\left\|x^{\mathbb{N}\backslash M}\right\|_{l_3}^3 + \sqrt{\|x^{\mathbb{N}\backslash M}\|_{l_3}^6 + 4\|x\|_{l_3}^3}\right)\left\|x^{\mathbb{N}\backslash M}\right\|_{l_3}^9 > 2\|x\|_{l_3}^6 - 2\|x\|_{l_3}^3\|x^M\|_{l_3}^6 \\
\Leftrightarrow\quad & \left(1 + \sqrt{1 + 4\frac{\|x\|_{l_3}^3}{\|x^{\mathbb{N}\backslash M}\|_{l_3}^6}}\right)\left\|x^{\mathbb{N}\backslash M}\right\|_{l_3}^{12} > 2\|x\|_{l_3}^6 - 2\|x\|_{l_3}^3\|x^M\|_{l_3}^6 \\
\Leftrightarrow\quad & 1 + \sqrt{1 + 4\frac{\|x\|_{l_3}^3}{\|x^{\mathbb{N}\backslash M}\|_{l_3}^6}} > 2\left(\frac{\|x\|_{l_3}^3}{\|x^{\mathbb{N}\backslash M}\|_{l_3}^6}\right)^2 - 2\frac{\|x\|_{l_3}^3}{\|x^{\mathbb{N}\backslash M}\|_{l_3}^6}.
\end{aligned} \tag{6.19}$$

Where $0 < \frac{\|x\|_{l_3}^3}{\|x^{\mathbb{N}\setminus M}\|_{l_3}^6} < \infty$ and $\|x\|_{l_3}^3 > 1$. Notice that

$$1 + \sqrt{1 + 4\frac{\|x\|_{l_3}^3}{\|x^{\mathbb{N}\setminus M}\|_{l_3}^6}} > 2.$$

So, if $\frac{\|x\|_{l_3}^3}{\|x^{\mathbb{N}\setminus M}\|_{l_3}^6}$ satisfies

$$0 < \frac{\|x\|_{l_3}^3}{\|x^{\mathbb{N}\setminus M}\|_{l_3}^6} \leq \frac{1+\sqrt{5}}{2},$$

then (6.19) holds. That is, if

$$\|x\|_{l_3} \leq \left(\frac{1+\sqrt{5}}{2}\right)^{\frac{1}{3}} \left\|x^{\mathbb{N}\setminus M}\right\|_{l_3}^2,$$

then,

$$\frac{\|x^M\|_{l_3}}{\|x\|_{l_3}} > \frac{1}{\|x^M\|_{l_3}\|a\|_{l_3}}. \qquad \square$$

Now, by Lemmas 6.13 and 6.14, we find the solutions of the generalized metric projection Π from $l_3$ onto its closed and convex cylinder $C_M(1)$. It is clear that, for any $x \in C_M(1)$, we have $\Pi_{C_M(1)}(x) = x$. So, we only consider some points in $l_3\setminus C_M(1)$.

**Theorem 6.15**. *For $x \in l_3\setminus C_M(1)$. Suppose x satisfies*

$$\|x\|_{l_3} \leq \left(\frac{1+\sqrt{5}}{2}\right)^{\frac{1}{3}} \left\|x^{\mathbb{N}\setminus M}\right\|_{l_3}^2, \tag{6.18}$$

*Let $b = b(x)$ and $a = a(x)$ be defined as in Lemma* 6.13. *Then*

$$\Pi_{C_M(1)}(x) = \frac{1}{\|x^M\|_{l_3}} x^M + bx^{\mathbb{N}\setminus M}.$$

*Proof*. Part (a) is clear. So, we only prove (b). For any given $x \in l_3\setminus C_M(1)$, $x$ satisfies $\|x\|_{l_3} \geq \|x^M\|_{l_3} > 1$. By condition (6.18), $x$ satisfies $\left\|x^{\mathbb{N}\setminus M}\right\|_{l_3} > 0$. Hence $b = b(x)$ and $a = a(x)$ are well defined as in Lemma 6.13. For any $y \in C_M(1)$ satisfying $\|y^M\|_{l_3} \leq 1$, by (6.16) and part (b) of Lemma 6.13 and Lemma 6.14, we have

$$\begin{aligned}
&\langle J(x) - J\left(\frac{1}{\|x^M\|_{l_3}} x^M + bx^{\mathbb{N}\setminus M}\right), \frac{1}{\|x^M\|_{l_3}} x^M + bx^{\mathbb{N}\setminus M} - y\rangle \\
&= \langle J(x) - J(a), \frac{1}{\|x^M\|_{l_3}} x^M + bx^{\mathbb{N}\setminus M} - y\rangle \\
&= \langle (J(x))^M - (J(a))^M, \frac{1}{\|x^M\|_{l_3}} x^M + bx^{\mathbb{N}\setminus M} - y\rangle \\
&= \left(\frac{\|x^M\|_{l_3}}{\|x\|_{l_3}} - \frac{1}{\|x^M\|_{l_3}\|a\|_{l_3}}\right)\langle J(x^M), \frac{1}{\|x^M\|_{l_3}} x^M - y^M\rangle \\
&= \left(\frac{\|x^M\|_{l_3}}{\|x\|_{l_3}} - \frac{1}{\|x^M\|_{l_3}\|a\|_{l_3}}\right)(\|x^M\|_{l_3} - \langle J(x^M), y^M\rangle) \\
&\geq \left(\frac{\|x^M\|_{l_3}}{\|x\|_{l_3}} - \frac{1}{\|x^M\|_{l_3}\|a\|_{l_3}}\right)\|x^M\|_{l_3}(1 - \|y^M\|_{l_3}) \\
&\geq 0, \text{ for all } y \in C_M(1).
\end{aligned} \tag{6.20}$$

Since

$$\left\| \left( \frac{1}{\|x^M\|_{l_3}} x^M \ + \ bx^{\mathbb{N}\backslash M} \right)^M \right\|_{l_3} = \left\| \frac{1}{\|x^M\|_{l_3}} x^M \right\|_{l_3} = 1,$$

this implies

$$\frac{1}{\|x^M\|_{l_3}} x^M \ + \ bx^{\mathbb{N}\backslash M} \in C_M(1).$$

By (6.20) and by using the basic variational principle of $\Pi_{C_M(1)}$, this implies that, for $x \in l_3\backslash C_M(1)$, if it satisfies (6.18), then

$$\Pi_{C_M(1)}(x) = \frac{1}{\|x^M\|_{l_3}} x^M + bx^{\mathbb{N}\backslash M}. \qquad \square$$

**Remarks 6.16**. From Theorem 6.15, we notice that the solutions of $\Pi_{C_M(1)}$ are very complicated, even for $p = 3$, $r = 1$ and under condition (6.18). We can imagine that, even in this special case, $\Pi'_{C_M(r)}$ is extremely complicated. So, we will omit the calculation for $\Pi'_{C_M(r)}$ in this paper.

### 6.5. The standard metric projection *P* onto closed and convex cylinders in $l_3$

In this subsection, we study the standard metric projection $P$ onto closed and convex cylinders in $l_p$, for $1 < p < \infty$. We calculate the values of the metric projection $P_{C_M(r)}$ from $l_p$ onto $C_M(r)$. In contrast with the results about the generalized metric projection $\Pi_{C_M(r)}$ studied in previous subsection, we will see that the solutions of the standard metric projection $P_{C_M(r)}$ onto closed and convex cylinders in $l_p$ are simpler than the solutions of $\Pi_{C_M(r)}$, which helps us to calculate the directional derivatives of $P_{C_M(r)}$.

**Theorem 6.17**. Let *p be a positive number with* $1 < p < \infty$. *For any* $x \in l_p$, *we have*

(a) $P_{C_M(r)}(x) = x$, *for* $x \in C_M(r)$;

(b) $P_{C_M(r)}(x) = \frac{r}{\|x^M\|_{l_p}} x^M + x^{\mathbb{N}\backslash M}$, *for* $x \in l_p\backslash C_M(r)$.

*Proof*. Part (a) is clear. So, we only prove (b). For any $x \in l_p\backslash C_M(r)$, $x$ satisfies $\|x\|_{l_p} \geq \|x^M\|_{l_p} > r$. Then, for any $y \in C_M(r)$ with $\|y^M\|_{l_r} \leq r$, by Lemma 6.1, we have

$$\begin{aligned}
&\langle J\left(x - \left(\frac{r}{\|x^M\|_{l_p}} x^M \ + \ x^{\mathbb{N}\backslash M}\right)\right), \frac{r}{\|x^M\|_{l_p}} x^M + x^{\mathbb{N}\backslash M} - y\rangle \\
&= \langle J\left(x^M + x^{\mathbb{N}\backslash M} - \left(\frac{r}{\|x^M\|_{l_p}} x^M \ + \ x^{\mathbb{N}\backslash M}\right)\right), \frac{r}{\|x^M\|_{l_p}} x^M + x^{\mathbb{N}\backslash M} - y\rangle \\
&= \langle J\left(x^M - \frac{r}{\|x^M\|_{l_p}} x^M\right), \frac{r}{\|x^M\|_{l_p}} x^M + x^{\mathbb{N}\backslash M} - (y^M + y^{\mathbb{N}\backslash M})\ \rangle \\
&= \left(1 - \frac{r}{\|x^M\|_{l_p}}\right)\langle J(x^M), \frac{r}{\|x^M\|_{l_p}} x^M - y^M\rangle \\
&= \left(1 - \frac{r}{\|x^M\|_{l_p}}\right)(r\|x^M\|_{l_p} - \langle J(x^M), y^M\rangle) \\
&\geq \left(1 - \frac{r}{\|x^M\|_{l_p}}\right)\|x^M\|_{l_p}(r - \ \|y^M\|_{l_p})
\end{aligned}$$

$\geq 0$, for all $y \in C_M(r)$. (6.21)

Since

$$\left\| \frac{r}{\|x^M\|_{l_p}} x^M \right\|_{l_p} = 1,$$

this implies

$$\frac{r}{\|x^M\|_{l_p}} x^M \in S_M(r) \subseteq C_M(r).$$

By (6.21) and by using the basic variational principle of $P_{C_M(r)}$, this implies that, for $x \in l_p \backslash C_M(r)$, we have

$$P_{C_M(r)}(x) = \frac{r}{\|x^M\|_{l_p}} x^M + x^{\mathbb{N}\backslash M}. \qquad \square$$

With the help of Theorem 6.17, in next theorem, we prove the Fréchet differentiability and calculate the directional derivatives of the metric projection $P_{C_M(r)}$ onto closed and convex cylinder $C_M(r)$ in $l_p$.

**Theorem 6.18**. *For any* $r > 0$, *we have,*

(a) $P_{C_M(r)}$ *is Fréchet differentiable on* $C_M^o(r)$ *such that, for any* $x \in C_M^o(r)$,

$$P'_{C_M(r)}(x)(h) = h, \text{ for any } h \in l_p \backslash \{\theta\};$$

(b) $P_{C_M(r)}$ *is Fréchet differentiable on* $l_p \backslash C_M(r)$ *such that, for any* $x \in l_p \backslash C_M(r)$,

$$P'_{C_M(r)}(x)(h) = \frac{r}{\|x^M\|_{l_p}} \left( h^M - \frac{\Psi(x^M,\ h^M)}{\|x^M\|_{l_p}} x^M \right) + h^{\mathbb{N}\backslash M}, \text{ for any } h \in l_p \backslash \{\theta\}.$$

*In particular,*

$$P'_{C_M(r)}(x)(x) = x^{\mathbb{N}\backslash M}, \text{ for any } x \in l_p \backslash C_M(r).$$

*Proof*. Proof of (a). Let $x \in C_M^o(r)$. Since $\|x^M\|_{l_p} < r$, then, for any $h \in l_p$ with $h \neq \theta$, there is $\delta > 0$, such that $\|(x + th)^M\|_{l_p} < r$, that is, $x + th \in C_M^o(r)$, for all $t \in (0, \delta)$, in which, $\delta$ can be chosen independently from $h$ with $\|h\|_{l_p} = 1$. By the explanations before Theorem 6.17, we have

$$\begin{aligned} & P'_{C_M(r)}(x)(h) \\ &= \lim_{t\downarrow 0} \frac{P_{C_M(r)}(x+th) - P_{C_M(r)}(x)}{t} \\ &= \lim_{t\downarrow 0, t<\delta} \frac{x+th-x}{t} \\ &= \lim_{t\downarrow 0, t<\delta} \frac{th}{t} \\ &= h. \end{aligned}$$

Since $\delta$ is chosen independently from $h$ with $\|h\|_{l_p} = 1$, this implies that the above limit is uniformly convergent with respect to $h$ satisfying $\|h\|_{l_p} = 1$. Then, it follows that the metric projection $P_{C_M(r)}$ is Fréchet differentiable at the point $x \in C_M^o(r)$.

Proof of (b). Let $x \in l_p \backslash C_M(r)$. Since $\|x^M\|_{l_p} > r$, we can choose a number $\lambda$ that is independent from $h$ with $\|h\|_{l_p} = 1$, such that, for any $t \in (0, \lambda)$ and any $h \in l_p$ with $\|h\|_{l_p} = 1$, we have

$$\|(x+th)^M\|_{l_p} > r, \text{ for any } t \in (0, \lambda) \text{ and for any } h \in l_p \text{ with } \|h\|_{l_p} = 1.$$

This implies

$$x + th \in l_p \backslash C_M(r), \text{ for any } t \in (0, \lambda) \text{ and for any } h \in l_p \text{ with } \|h\|_{l_p} = 1.$$

By Theorem 6.17 and Lemma 6.1, we have

$$\begin{aligned}
&P'_{C_M(r)}(x)(h)\\
&= \lim_{t\downarrow 0} \frac{P_{C_M(r)}(x+th) - P_{C_M(r)}(x)}{t}\\
&= \lim_{t\downarrow 0, t<\lambda} \frac{\left(\frac{r}{\|(x+th)^M\|_{l_p}}(x+th)^M + (x+th)^{\mathbb{N}\backslash M}\right) - \left(\frac{r}{\|x^M\|_{l_p}}x^M + x^{\mathbb{N}\backslash M}\right)}{t}\\
&= \lim_{t\downarrow 0, t<\lambda} \frac{\frac{rth^M}{\|(x+th)^M\|_{l_p}} + \left(\frac{r}{\|(x+th)^M\|_{l_p}} - \frac{r}{\|x^M\|_{l_p}}\right)x^M + (th)^{\mathbb{N}\backslash M}}{t}\\
&= \lim_{t\downarrow 0, t<\lambda} \frac{rh^M}{\|(x+th)^M\|_{l_p}} + h^{\mathbb{N}\backslash M} + \lim_{t\downarrow 0, t<\lambda} \frac{\left(\frac{r}{\|(x+th)^M\|_{l_p}} - \frac{r}{\|x^M\|_{l_p}}\right)x^M}{t}\\
&= \frac{r}{\|x^M\|_{l_p}} h^M + h^{\mathbb{N}\backslash M} - \frac{r}{\|x_M\|^2}\Psi(x^M, \ h^M)x^M\\
&= \frac{r}{\|x^M\|_{l_p}}\left(h^M - \frac{\Psi(x^M, \ h^M)}{\|x^M\|_{l_p}}x^M\right) + h^{\mathbb{N}\backslash M}. \qquad \square
\end{aligned}$$